\documentclass[reqno, 11pt, letterpaper ]{amsart}

\usepackage{amsmath}
\usepackage{amsfonts}
\usepackage{amssymb}
\usepackage{graphicx}
\usepackage{amsthm,color,yfonts,cite}
\usepackage{paralist}
\usepackage{hyperref}
\usepackage{physics}
\usepackage{bbm}
 \usepackage{bm}
\usepackage{comment}
\usepackage{amsaddr}
\usepackage{mathtools}

\usepackage{tikz}
\usepackage{tikz-cd}
\usetikzlibrary{cd}
\usetikzlibrary{arrows.meta}
\usetikzlibrary{matrix, arrows, positioning}
\usetikzlibrary{shadows}

\newtheorem{theorem}{Theorem}
\newtheorem{definition}[theorem]{Definition}
\newtheorem{proposition}[theorem]{Proposition}
\newtheorem{lemma}[theorem]{Lemma}
\newtheorem{corollary}[theorem]{Corollary}

\theoremstyle{remark}
\newtheorem{remark}[theorem]{Remark}

\makeatletter
\newcommand*{\rom}[1]{\expandafter\@slowromancap\romannumeral #1@}
\makeatother

\newcommand{\ls}{\lesssim}

\newcommand{\R}{\mathbb{R}}

\newcommand{\ep}{\epsilon}
\newcommand{\al}{\alpha}

\newcommand{\CB}{\mathcal{B}}

\newcommand{\CR}{\mathcal{R}}

\newcommand{\CU}{\mathcal{U}}

\newcommand{\CX}{\mathcal{X}}
\newcommand{\CY}{\mathcal{Y}}
\newcommand{\CZ}{\mathcal{Z}}

\newcommand{\FL}{\mathfrak{L}}

\newcommand{\FS}{\mathfrak{S}}

\newcommand{\ga}{\gamma}

\newcommand{\ps}{\psi}
\newcommand{\ph}{\varphi}

\renewcommand{\th}{\theta}

\newcommand{\si}{\sigma}

\newcommand{\et}{\eta_0}

\newcommand{\rh}{\rho}

\newcommand{\pl}{\partial}
\newcommand{\wt}{\widetilde}

\newcommand{\Ck}[1]{\left\{#1\right\}}

\newcommand{\K}[1]{\left(#1\right)}

\newcommand{\No}[1]{\left\| #1 \right\|}
\newcommand{\I}{\infty}
\newcommand{\sd}{\langle \hb \nabla \rangle}
\newcommand{\bx}{\langle x \rangle}

\newcommand{\tw}{\frac{1}{2}}
\newcommand{\na}{\nabla}
\newcommand{\hb}{\hbar}

\newcommand{\rg}{\rangle}
\renewcommand{\lg}{\langle}
\newcommand{\bt}{\lg t \rg}

\newcommand{\J}{\mathcal{J}}
\newcommand{\U}{\CU^\hb(t)}
\newcommand{\Uinv}{\CU^\hb(t)^*}
\newcommand{\M}{e^{\frac{i|x|^2}{2\hb t}}}
\newcommand{\Minv}{e^{-\frac{i|x|^2}{2\hb t}}}

\newcommand{\Weyl}{\mathbf{Weyl}^\hbar}
\newcommand{\Wig}{\mathbf{Wig}^\hbar}

\newcommand{\Tp}{\mathbf{Toep}^\hb}
\newcommand{\Hs}{\mathbf{Hus}^\hb}
\newcommand{\E}{\CU(t)}

\usepackage{wrapfig}
\usepackage{tikz}
\usetikzlibrary{arrows,calc,decorations.pathreplacing}
\definecolor{light-gray1}{gray}{0.90}
\definecolor{light-gray2}{gray}{0.80}
\definecolor{light-gray3}{gray}{0.60}

\numberwithin{equation}{section}

\numberwithin{theorem}{section}

\numberwithin{table}{section}

\numberwithin{figure}{section}

\ifx\pdfoutput\undefined
  \DeclareGraphicsExtensions{.pstex, .eps}
\else
  \ifx\pdfoutput\relax
    \DeclareGraphicsExtensions{.pstex, .eps}
  \else
    \ifnum\pdfoutput>0
      \DeclareGraphicsExtensions{.pdf}
    \else
      \DeclareGraphicsExtensions{.pstex, .eps}
    \fi
  \fi
\fi

\title[Semi-classical limit of quantum scattering states]{Semi-classical limit of quantum scattering states for the nonlinear Hartree equation}

\author[S. Hadama]{Sonae Hadama}
\address{The University of Osaka, 1-1 Machikaneyama-cho, Toyonaka-shi, Osaka 560-0043, Japan}
\email{hadama.sonae.sci@osaka-u.ac.jp}

\author[Y. Hong]{Younghun Hong}
\address{Department of Mathematics, Chung-Ang University, Seoul 06974, Korea}
\email{yhhong@cau.ac.kr (Corresponding Author)}

\date{\today}
\begin{document}

\begin{abstract}
This article concerns the long-time dynamics of quantum particles in the semi-classical regime. First, we show that for the nonlinear Hartree equation with short-range interaction potential, small-data solutions obey dispersion bounds (Theorem \ref{thm: uniform bounds for nonlinear solutions}) and they scatter (Corollary \ref{cor: uniform small-data scattering}), where the smallness conditions and the bounds are independent of the small parameter $\hbar\in(0,1]$ representing the reduced Planck constant. Then, taking the semi-classical limit $\hbar\to0$, we prove that the Wigner transforms of such quantum scattering states converge weakly to the corresponding classical scattering states for the Vlasov equation (Theorem \ref{thm: semiclassical limit}). As a direct consequence, we establish small-data scattering for the Vlasov equation without assuming regularity on initial data (Theorem \ref{thm:scattering Vlasov}). Our analysis is based on a new uniform dispersion estimate for  the free Schr\"odinger flow (Proposition \ref{prop:free dispersive}), which is simple but crucial for including singular interaction potentials such as inverse power-law potential $\frac{1}{|x|^a}$ with $1<a<\frac{5}{3}$.
\end{abstract}

\maketitle

\section{Introduction}

The nonlinear Hartree equation and the Vlasov equation are two fundamental models describing the mean-field dynamics of a large number of quantum and classical particles, respectively. By the quantum-classical correspondence, the Vlasov equation is obtained as the semi-classical limit of the nonlinear Hartree equation. On the other hand, for both equations, it is known that if the pair interactions are short-range, small-data solutions behave asymptotically linearly; in other words, they scatter.
The purpose of this paper is to examine the consistency of these long-time dynamics between the quantum and the classical mean-field models, which is more than convergence on a finite time interval.

\subsection{Nonlinear Hartree equation}

In the Heisenberg picture, the nonlinear Hartree equation is given by
\begin{equation}\label{eq: NLH}
i\hbar\partial_t\gamma^\hbar=\bigg[-\frac{\hbar^2}{2}\Delta+\Phi^\hbar,\gamma^\hbar\bigg],
\end{equation}
where $\gamma^\hbar=\gamma^\hbar(t)$ is a time-dependent bounded self-adjoint operator acting on $L^2(\mathbb{R}^3)$, $[A, B]=AB-BA$ is the Lie bracket, and $\hbar>0$ represents the reduced Planck constant. With a slight abuse of notation, $A(x,x')$ denotes the kernel of an integral operator $A$, that is, 
$$(A\phi)(x)=\int_{\mathbb{R}^3}A(x,x')\phi(x')dx',$$
and we formally\footnote{Rigorously, $\rho_A^\hbar=(2\pi\hbar)^3A(x,x)$ is well-defined for smooth finite-rank operators. Then, it can be extended to a larger class of operators as a natural extension by the linear estimates.} define the (semi-classically rescaled) density function for an operator $A$ by 
$$\rho^\hbar_A:=(2\pi\hbar)^3A(x,x).$$
Then, the self-generated mean-field potential in \eqref{eq: NLH} is given by
\begin{equation}\label{eq: quantum mean-field potential}
\Phi^\hbar:=w \ast \rh_{\gamma^\hbar}^\hbar,
\end{equation}
where $w:\mathbb{R}^3\to\mathbb{R}$ is the pair interaction potential among particles.

The equation \eqref{eq: NLH} is a mean-field approximation model for quantum interacting particles, and it can be used for both bosons and fermions. Indeed, if quantum observables are restricted to finite-rank projectors, i.e., $\gamma^\hbar=\sum_{j=1}^N |\phi_j^\hbar\rangle\langle\phi_j^\hbar|$ with an orthonormal set $\{\phi_j^\hbar\}_{j=1}^N\subset L^2(\mathbb{R}^3)$, then \eqref{eq: NLH} is equivalent to a system of coupled equations 
\begin{equation}\label{eq: NLH, coupled}
i\hbar\partial_t\phi_j^\hbar=-\frac{\hbar^2}{2}\Delta\phi_j^\hbar+\bigg((2\pi\hbar)^3\sum_{k=1}^N w*|\phi_k^\hbar|^2\bigg)\phi_j^\hbar. 
\end{equation}
In particular, when $N=1$, this is a well-known mean-field model for Bose-Einstein condensates (see \cite{Erdos et al 2006, Erdos et al 2007, Erdos et al 2009, Erdos et al 2010} for a rigorous derivation). For fermions, the Hartree-Fock equation is derived from the many-body linear Schr\"odinger equation (see, for example, \cite{Benedikter et al 2014, Benedikter et al 2016 derivation}). In the semi-classical regime, it is often reduced further to \eqref{eq: NLH} dropping the exchange term, while an additional constraint $0\leq \gamma^\hbar\leq 1$ is imposed due to Pauli's exclusion principle.

Even though the two formulations \eqref{eq: NLH} and \eqref{eq: NLH, coupled} share many common analytic features, in this paper, we focus on the operator-valued equation \eqref{eq: NLH}, since it has several significant advantages.  Indeed, the equation \eqref{eq: NLH} is much more natural in connection with the long-time dynamics for the Vlasov equation; the finite system \eqref{eq: NLH, coupled} with fixed $N\geq1$ yields only a finite-$N$ sum of non-dispersive Dirac measures in the semi-classical limit (see \cite{Lions Paul 1993}). In addition, the operator-valued equation \eqref{eq: NLH} admits some physically important quantum states having infinite particles, and much richer dynamics can be described around such infinite-particle states (see \cite{Lewin Sabin 2015, Lewin Sabin 2014, Chen et al 2017, Chen et al 2018, Hadama 2023, Hadama Yamamoto 2023, Hadama Hong 2024, You 2024, Smith 2024}).

The global well-posedness of \eqref{eq: NLH} is shown in \cite{Bove et al 1974, Bove et al 1976, Chadam 1976, Zagatti 1992}. For long-time dynamics with fixed $\hbar>0$, a small-data scattering result is given in \cite{Pusateri Sigal 2021} for short-range interactions (see \eqref{eq: quantum scattering}). It is shown in this article that for \eqref{eq: NLH}, small-data solutions scatter (Corollary \ref{cor: uniform small-data scattering}). Recently, a modified scattering for long-range interactions is proved  \cite{Nguyen You 2024}.

\subsection{Vlasov equation}
The Vlasov equation is a kinetic model that describes the mean-field dynamics of classical collisionless particles. It is given by \begin{equation}\label{eq: Vlasov}
\partial_t f+p\cdot\nabla_q f-\nabla \Phi\cdot\nabla_pf=0,
\end{equation}
where $f=f(t,q,p): I(\subset\mathbb{R})\times\mathbb{R}^3\times\mathbb{R}^3\to[0,\infty)$ is a distribution function on the phase space and 
\begin{equation}\label{eq: classical mean-field potential}
\Phi=w*\rho_f
\end{equation}
is the mean-field potential, that is, the convolution of the interaction potential $w$ and the density function
$$\rho_f:=\int_{\mathbb{R}^3}f(\cdot,p)dp.$$
In the literature, due to its physical importance, the Vlasov-Poisson equation, i.e., \eqref{eq: Vlasov} with Coulomb interaction $w=\pm\frac{1}{|x|}$, has drawn special attention.  For the initial-value problem, global solutions are constructed under mild conditions on initial data \cite{Lions Perthame 1991, Pfaffelmoser 1992}. On the other hand, in Bardos and Degond \cite{Bardos Degond 1985}, by a different approach, global well-posedness of small-data solutions, together with a time decay bound for densities, has been established. Then, the modified scattering is proved in Choi-Kwon \cite{Choi Kwon 2016} (see also Ionescu-Pausader-Wang-Widmayer \cite{IPWW 2022}).

For the equation \eqref{eq: Vlasov} with a general potential $w$, the approach in Bardos-Degond can be employed for the long-time dynamics. Similarly to the quantum analogue, it is shown in \cite{Choi Ha 2011} that if $w$ decays faster than the Coulomb potential, then small-data solutions decay in time and scatter (see \eqref{eq: definition of kinetic scattering}), while non-scattering is shown in the long-range interaction case.

\subsection{Semi-classical limit}

In mathematical physics, one of the central topics is the rigorous derivation of a mesoscopic kinetic model from a microscopic quantum mechanical model by semi-classical analysis (and that of a macroscopic fluid model by the hydrodynamic limit). In this article, we are particularly concerned with the semi-classical limit from the nonlinear Hartree equation to the Vlasov equation. To the best of the authors' knowledge, the first rigorous proof was given independently by Lions and Paul \cite{Lions Paul 1993} and Markowich and Mauser \cite{Markowich Mauser 1993}, where the Vlasov equation is derived from the Hartree equation in a weak sense, including singular interactions, but a robust mathematical framework for semi-classical analysis was also introduced. Since then, the semi-classical limits in the weak convergence regime were developed by \cite{Gasser et al 1998, Figalli et al 2012}.

In recent years, there has been significant progress on this problem. In particular, the weak convergence to the Vlasov equation is upgraded to strong convergence with an explicit rate of convergence when the interaction potential is regular enough \cite{Pezzotti Pulvirenti 2009, Athanassoulis et al 2011, Amour et al 2013a, Amour et al 2013b, Benedikter et al 2016 semiclassical, Golse Paul 2017}. Then, the theory has been developed to include singular interaction cases \cite{Chong et al 2022, Chong et al 2023, Chong et al 2024, Lafleche 2019, Lafleche 2021, Lafleche Saffirio 2023, Saffirio 2020a, Saffirio 2020b}. Most notably, in the recent work of Lafleche and Saffirio \cite{Lafleche Saffirio 2023} and Chong, Lafleche and Saffirio \cite{Chong et al 2023}, an affirmative answer was given for the strong convergence to the Vlasov-Poisson equation in the semi-classical limit. In addition, the Vlasov-Poisson equation is derived from the quantum many-body problem via the mean-field limit \cite{Chong et al 2024}. 
In the setting of the quantum analogue of the Landau damping, we refer to \cite{Lewin Sabin 2020, Smith 2024}.

 

\subsection{Basic notations}\label{sec: notations}

In this paper, we consider the semi-classical limit $\hbar\to0$. Therefore, it is important to employ the semi-classically rescaled definitions and quantities for quantum observables. 

\subsubsection{Semi-classically scaled operator norms}
Let $\CB(L^2(\mathbb{R}^3))$ be the space of all bounded linear operators on $L^2(\R^3)$.
For $r \in [1,\I)$, the Schatten-$r$ class $\FS^r$ is defined as the Banach subspace of $\CB(L^2(\mathbb{R}^3))$ equipped with the norm
$$\|A\|_{\FS^r}=(\Tr(|A|^r))^{\frac{1}{r}},$$
where $|A|:= \sqrt{A^* A}$. However, the semi-classically rescaled norm given by
$$\|A\|_{\mathfrak{L}_\hbar^r}:=(2\pi\hbar)^{\frac{3}{r}}\|A\|_{\mathfrak{S}^r}$$
is convenient to use. Indeed, for fixed $\hbar>0$, the $\mathfrak{L}_\hbar^r$-norm is equivalent to the Schatten $r$-norm, but in the semi-classical regime, formally, it corresponds to the $L_{q,p}^r(\mathbb{R}^6)$ norm on the phase space $\mathbb{R}^6=\mathbb{R}^3\times\mathbb{R}^3(\ni (q,p))$ via the Wigner transform (see \eqref{eq: Wigner transform}). Accordingly, it is natural to define the quantum weighted (in both position and momentum) norm by 
$$\|\gamma\|_{\mathfrak{L}_\hbar^{r,\sigma}}:=(2\pi\hbar)^{\frac{3}{r}}\|\langle \hbar\nabla\rangle^{\sigma}\gamma\langle \hbar\nabla\rangle^{\sigma}\|_{\mathfrak{S}^r}+(2\pi\hbar)^{\frac{3}{r}}\|\langle x\rangle^{\sigma}\gamma\langle x\rangle^{\sigma}\|_{\mathfrak{S}^r},$$
which again, by the Wigner transform, corresponds to the weighted norm
$$\|f\|_{L_{q,p}^{r,\sigma}(\mathbb{R}^{6})}:=\|\langle q\rangle^\sigma f\|_{L_{q,p}^{r}(\mathbb{R}^{6})}+\|\langle p\rangle^\sigma f\|_{L_{q,p}^{r}(\mathbb{R}^{6})}.$$
For more details, we refer to Lafleche \cite{Lafleche 2024}.

\subsubsection{Linear propagators}\label{sec: linear propagators}
We introduce the notations for the three key linear propagators. Given a time-dependent potential $V(t)$, let $\mathcal{U}_V^\hbar(t,t_0)$ be the linear flow associated with the linear Schr\"odinger equation perturbed by the potential $V$. Precisely, 
$$u(t)=\mathcal{U}_V^\hbar(t,t_0)\varphi$$
is the solution to 
\begin{equation}\label{eq: linear Schrodinger}
i\hbar\partial_t u=\bigg(-\frac{\hbar^2}{2}\Delta +V(t)\bigg)u
\end{equation}
with $u(t_0)=\varphi$. In the free case $V\equiv0$, we denote
$$\CU^\hb(t):= e^{\frac{it\hbar}{2}\Delta}.$$ Finally, we let
\begin{equation}\label{eq: free transport flow, definition}
(\E f_0)(q,p):=f_0(q-tp,p)
\end{equation}
denote the solution to the free transport equation
$$\pl_t f + p \cdot \na_q f=0$$
with initial data $f(0,q,p)=f_0(q,p)$.

\subsubsection{Vector field}\label{sec: vector-field}
In our analysis, the vector field
$$\J_t^\hbar:=x+it\hbar\nabla,$$
associated with the semi-classically scaled free Schr\"odinger equation, will be frequently employed in that it arises naturally in studying scattering problems in a weighted space. For its properties, we refer to \cite{Ginibre et al 1994, Hayashi Naumkin 1998, Hayashi et al 1998} and references therein. Among many others, the identity 
\begin{align}\label{J to derivative}
	\J_t^\hb = \U x \Uinv = \M (it\hb\na) \Minv
\end{align}
is particularly useful. From this identity, the fractional power of the vector field can be defined as 
\begin{equation}\label{eq: factional power of J}
	\begin{aligned}
		&|\J_t^\hbar|^\sigma := \U |x|^\si \Uinv = e^{\frac{i|x|^2}{2t\hbar}}|t\hbar\nabla|^\sigma e^{-\frac{i|x|^2}{2t\hbar}} \\
		&\langle \J_t^\hbar\rangle^\sigma := \U \bx^\sigma \Uinv =e^{\frac{i|x|^2}{2t\hbar}}\langle t\hbar\nabla\rangle^\sigma e^{-\frac{i|x|^2}{2t\hbar}}
	\end{aligned}
\end{equation}

\section{Main results}

Throughout this article, we assume that the pair interaction $w$ for both the quantum and classical models \eqref{eq: NLH} and \eqref{eq: Vlasov} satisfies the following inverse-power law:
\begin{equation}\label{eq: potential assumption}
w(x)=\pm\frac{1}{|x|^a}\quad\textup{with }1<a<\frac{5}{3}.
\end{equation}
\begin{remark}
$(i)$ The Vlasov equation \eqref{eq: Vlasov} with the interaction potential $w=\pm\frac{1}{|x|^a}$ is called the Vlasov-Riesz system. In particular, it is the Vlasov-Poisson system when $a=1$, that is, the most important Coulomb interaction case. \\
$(ii)$ The long-range case $a\leq 1$ is excluded, because there is no scattering. In the important borderline case $a=1$, the modified scattering is established for \eqref{eq: NLH} in \cite{Nguyen You 2024}, and for \eqref{eq: Vlasov} in \cite{Choi Kwon 2016}. It would also be interesting to investigate their semi-classical correspondence, but it will be postponed to our future work. The more singular interaction case $a\geq \frac{5}{3}$ is not included in this paper due to technical difficulties (see the proof of Lemma \ref{lemma: interpolation inequalities for potentials}).\\
$(iii)$ $w$ is restricted to the inverse power-law potential to avoid non-essential technical details. Our analysis can be easily extended to general potentials of the form $w=w_1+w_2$, where $w_j\in L^{\frac{3}{a_j},\infty}(\mathbb{R}^3)$ and $\nabla w_j\in L^{\frac{3}{a_j+1},\infty}(\mathbb{R}^3)$ for some $1<a_j<\frac{5}{3}$, because only proper upper bounds are required for potentials.\\
$(iv)$ We also fix the dimension to $d=3$, which corresponds to the physically most relevant case. This restriction is motivated not only by numerical simplicity but also to allow a comparison with the endpoint model $a=1$. In this case, strong convergence toward the Vlasov-Poisson system has been established on compact time intervals \cite{Lafleche Saffirio 2023}, and modified scattering is known \cite{Choi Kwon 2016, IPWW 2022}.
\end{remark}

\subsection{Uniform dispersion estimates for the nonlinear Hartree equation}\label{sec: uniform bounds}

To compare the quantum and classical long-time dynamics in the semi-classical regime, it is crucial to capture dispersive effects for the NLH \eqref{eq: NLH} that can be measured independently of $\hbar\in(0,1]$. Our first main result provides a new uniform global-in-time dispersion estimate and uniform norm bounds for small-data solutions to the nonlinear Hartree equation \eqref{eq: NLH}, which we think has its own interest from an analysis perspective. 

\begin{theorem}[Uniform global-in-time bounds for the nonlinear Hartree equation]\label{thm: uniform bounds for nonlinear solutions}
If $w$ satisfies \eqref{eq: potential assumption} for some $1<a<\frac{5}{3}$, then for any sufficiently small $\epsilon>0$, there exists $\eta_0=\eta_0(a,\epsilon)>0$ such that the following holds. We assume that 
\begin{equation}\label{eq: initial data assumption}
\gamma^\hbar_0\geq0\quad\textup{and}\quad\sup_{\hbar\in(0,1]}\Big\{\|\gamma^\hbar_0\|_{\mathfrak{L}_\hbar^1}+\|\gamma^\hbar_0\|_{\mathfrak{L}_\hbar^{\frac{3}{2-a-\epsilon},\frac{1+a+2\epsilon}{2}}}\Big\}\leq\eta_0.
\end{equation}
Then, the solution $\gamma^\hbar(t)\in C(\mathbb{R}; \mathfrak{S}^1)$ to the nonlinear Hartree equation \eqref{eq: NLH} with initial data $\gamma_0^\hbar$ exists globally in time, and its density and potential function satisfy the decay bounds
\begin{equation}\label{eq: dispersion estimate}
\sup_{t \in \R} \Ck{\|\rh^\hb_{\ga^\hb}(t)\|_{L^1(\R^3)}+\bt^{1+a+\epsilon}\|\rho^\hbar_{\gamma^\hbar}(t)\|_{L^{\frac{3}{2-a-\epsilon}}(\mathbb{R}^3)}}\lesssim \eta_0,
\end{equation}
\begin{equation}\label{eq: dispersion estimate'}
\sup_{t \in \R} \bt^{1+a}\|\nabla\Phi^\hbar(t)\|_{L^\infty(\mathbb{R}^3)}\lesssim \eta_0,
\end{equation}
as well as the global-in-time norm bounds
\begin{equation}\label{eq: weighted norm uniform bound}
\begin{aligned}
&\big\|\langle\J_t^\hb\rangle^{\frac{1+a+2\epsilon}{2}}\gamma^\hbar(t) \langle \J_t^\hb\rangle^{\frac{1+a+2\epsilon}{2}}\big\|_{C_t(\mathbb{R}; \mathfrak{L}_{\hbar}^{\frac{3}{2-a-\epsilon}})}\\
&\qquad \qquad \qquad +\big\|\langle \hbar\nabla\rangle^{\frac{1+a+2\epsilon}{2}}  \gamma^\hbar(t)\langle \hbar\nabla\rangle^{\frac{1+a+2\epsilon}{2}}\big\|_{C_t(\mathbb{R}; \mathfrak{L}_\hbar^{\frac{3}{2-a-\epsilon}})}\lesssim\eta_0,
\end{aligned}
\end{equation}
where the implicit constants in \eqref{eq: dispersion estimate}, \eqref{eq: dispersion estimate'} and \eqref{eq: weighted norm uniform bound} do not depend on $\hbar\in(0,1]$.
\end{theorem}

\begin{remark}
$(i)$ Recently, Smith proved uniform-in-$\hbar$ decay estimates for the nonlinear Hartree equation in the Landau damping setting \cite{Smith 2024}, which might include our result when $w$ is sufficiently regular. On the other hand, our main result includes singular interaction potentials.\\
$(ii)$ Theorem \ref{thm: uniform bounds for nonlinear solutions} includes the short-range interaction case $a>1$. The constraint $a<\frac{5}{3}$ is due to technical difficulties\footnote{In \eqref{eq: why a<5/3}, we need $\frac{1+3a+2\epsilon}{2}<3$}. Indeed, for larger $a\geq\frac{5}{3}$, it seems required to obtain a uniform decay bound for derivatives of density functions, but currently we are unable to prove this.
\end{remark}

Moreover, assuming more on the initial data, we prove stronger dispersion estimates.

\begin{proposition}[Improved dispersion]\label{prop: upgraded uniform bounds}
Under the notation and assumptions in Theorem \ref{thm: uniform bounds for nonlinear solutions}, the following hold.
\begin{enumerate}[$(i)$]
\item If $1<a<\frac{3}{2}$ and if we further assume 
$\displaystyle\sup_{\hbar\in(0,1]}\|\gamma^\hbar_0\|_{\mathfrak{L}_\hbar^{\infty,\frac{3}{2}+\epsilon}}\leq\eta_0$, then 
$$\|\rho^\hbar_{\gamma^\hbar}(t)\|_{L^\infty(\mathbb{R}^3)}\lesssim\frac{\eta_0}{\langle t\rangle^3}.$$
\item
If $\frac{3}{2}\leq a<\frac{5}{3}$ and if we further assume 
$\displaystyle\sup_{\hbar\in(0,1]}\|\gamma^\hbar_0\|_{\mathfrak{L}_\hbar^{\frac{3}{2a-3+2\epsilon},\frac{6-2a-\epsilon}{2}}}\leq\eta_0$, then 
$$\|\rho^\hbar_{\gamma^\hbar}(t)\|_{L^{\frac{3}{2a-3+2\ep}}(\R^3)}
\lesssim\frac{\eta_0}{\langle t\rangle^{6-2a-2\ep}}.$$
\end{enumerate}
\end{proposition}

For the proof, it is important to note that Theorem \ref{thm: uniform bounds for nonlinear solutions} cannot be obtained by the standard perturbative argument involving the Strichartz estimates for the linear Schr\"odinger equation $i\partial_t u=-\Delta u$ (see Cazenave \cite{Cazenave}, for instance). Indeed, it is easy to see from the Duhamel formula
\begin{equation}\label{eq: crude Duhamel formula}
\gamma^\hbar(t)=\U\gamma^\hbar(0)\Uinv -\frac{i}{\hbar}\int_0^t \CU^\hb(t-t') [\Phi^\hbar,\gamma^\hbar](t') \CU^\hb(t-t')^*  dt'
\end{equation}
that the dispersion from the propagator $e^{\frac{it\hbar}{2}\Delta}$ becomes weaker for smaller $\hbar>0$. More importantly, the integral term cannot be considered as a perturbation of the linear part due to the large $\frac{1}{\hbar}$ factor. To overcome this problem, we introduce a simple dispersion estimate for the free Schr\"odinger flow (Proposition \ref{prop:free dispersive}), but we also find a way to implement it into the wave operator formulation of the nonlinear Hartree equation \eqref{eq: linear propagator form NLH} to resolve the large $\frac{1}{\hbar}$-factor issue in \eqref{eq: crude Duhamel formula} (Section \ref{sec: uniform bounds for NLH}).

As a direct application of Theorem \ref{thm: uniform bounds for nonlinear solutions}, we prove small-data scattering. 

\begin{corollary}[Small-data scattering]\label{cor: uniform small-data scattering}
A small-data global solution $\gamma^\hbar(t)\in C(\mathbb{R}; \mathfrak{L}_\hbar^1)$ to the NLH \eqref{eq: NLH}, constructed in Theorem \ref{thm: uniform bounds for nonlinear solutions}, scatters in $\mathfrak{L}_\hbar^1$ as $t\to\pm\infty$ in the following sense: There exists $\ga^\hb_\pm \in \FS^1$ such that 
\begin{equation}\label{eq: quantum scattering}
    \|\CU^\hb(-t) \ga^\hb(t) \CU^\hb(t) -\ga^\hb_\pm\|_{\mathfrak{L}_\hbar^1}\ls \frac{\eta_0^2}{\hb\bt^{a-1}} \to 0 \text{ as } t\to\pm \I.
\end{equation}
\end{corollary}

\begin{remark}
$(i)$ It is important to note that in Corollary \ref{cor: uniform small-data scattering}, smallness condition does not depend on $\hbar\in(0,1]$. Indeed, if one tries to prove small-data scattering for the NLH \eqref{eq: NLH} just by following standard argument in Cazenave \cite[Chapter 7]{Cazenave}, one would see that the smallness condition may depend on $\hbar>0$. The uniform decay estimate (Theorem \ref{thm: uniform bounds for nonlinear solutions}) is crucially used to resolve this issue.\\
$(ii)$ In our proof, the convergence to the scattering states as $t \to \infty$ in \eqref{eq: quantum scattering} depends on $\hbar \in (0,1]$. While we expect that \eqref{eq: quantum scattering} holds uniformly for $\hbar \in (0,1]$, we are not able to prove this at present. We leave this issue for future work.
\end{remark}

\subsection{Semi-classical limit of quantum scattering states}

Next, we are concerned with the semi-classical limit of small-data scattering solutions $\gamma^\hbar(t)$, constructed in Theorem \ref{thm: uniform bounds for nonlinear solutions}, toward those to the Vlasov equation defined as follows.
\begin{definition}[Global scattering solution to the Vlasov equation]\label{definition: scattering solutions to Vlasov}
A distribution function $f(t)=f(t,q,p)\geq0$ is called an $H^{-1}$-global solution to the Vlasov equation \eqref{eq: Vlasov} with initial data $f_0$ if $g(t)=\mathcal{U}(-t)f(t)$ obeys the Vlasov equation in a moving reference frame 
\begin{equation}\label{eq:meaning of solution}
\begin{aligned}
g(t)&=f_0 + \int_0^t  \CU(-t_1) \Big\{(\na_q w \ast \rh_f) \cdot \nabla_pf\Big\}(t_1) dt_1\in C(\mathbb{R}_t;H_{q,p}^{-1}(\mathbb{R}^6))\\
&=f_0+\int_0^t (-t_1\na_q+\na_p)\CU(-t_1) \Big\{\nabla_q (w*\rho_f)\cdot f\Big\}(t_1)dt_1,
\end{aligned}
\end{equation}
where $\mathcal{U}(t)$ is the free transport flow (see \eqref{eq: free transport flow, definition}). It is said to scatter in $H_{q,p}^{-1}(\mathbb{R}^6)$ if there exist $f_\pm \in H^{-1}_{q,p}(\R^6)$ such that
\begin{equation}\label{eq: definition of kinetic scattering}
\lim_{t\to\pm\infty}\|\mathcal{U}(-t)f(t)-f_\pm\|_{H^{-1}_{q,p}(\R^6)}=\lim_{t\to\pm\infty}\|g(t)-f_\pm\|_{H^{-1}_{q,p}(\R^6)}=0.
\end{equation}
\end{definition}

In order to compare quantum observables with functions on the phase space $\mathbb{R}^3\times\mathbb{R}^3$, we employ the Wigner transform, a well-known fundamental tool in semi-classical analysis (see \cite{Lions Paul 1993} and Section \ref{sec: Wigner transform and Weyl quantization} for more details).

\begin{definition}[Wigner transform]\label{definition: Wigner transform}
For an operator $\gamma\in \mathfrak{L}_\hbar^2$ with integral kernel $\gamma(x,x')$, its Wigner transform is given by  
$$\mathbf{Wig}^\hbar[\gamma](q,p):=\int_{\mathbb{R}^3}\gamma\bigg(q+\frac{y}{2},q-\frac{y}{2}\bigg) e^{-\frac{iy\cdot p}{\hbar}}dy.$$
\end{definition}

Finally, we recall the weak topology of the separable space $L_{q,p}^2(\mathbb{R}^6)$ (see Brezis\cite[Chapter 3]{Brezis}). 

\begin{definition}[Metric for the weak topology]\label{definition: weak topology}
For $R>0$, a ball $B_R := \{f \in L^2_{q,p}(\mathbb{R}^6) : \|f\|_{L^2_{q,p}} \le R\}$ with the weak topology is metrizable by the metric 
\begin{equation}\label{eq: weak-topology metric}
\mathbf{d}^{w}(g_1,g_2):=\sum_{n=1}^\infty 2^{-n}\big|\langle \psi_n|g_1-g_2\rangle_{L_{q,p}^2(\mathbb{R}^6)}\big|
\end{equation}
for all $g_1,g_2\in B_R$, where $\{\psi_n\}_{n=1}^\infty$ is a dense subset of the unit ball in $L_{q,p}^2(\mathbb{R}^6)$. This means that $g_j\rightharpoonup g$ in $L_{q,p}^2(\mathbb{R}^6)$ if and only if $\mathbf{d}^w(g_j,g)\to0$ as $j\to\infty$. Therefore, the space $C_t(I; w-B_R)$ can be considered as the collection of $B_R$-valued functions which are continuous with respect to the metric $\mathbf{d}^w$, equipped with the metric 
$$\mathbf{d}_I^w\big(g_1(t),g_2(t)\big):=\sup_{t\in I}\mathbf{d}^w(g_1(t),g_2(t)).$$
\end{definition}

Our second main result shows the correspondence between quantum and classical scattering states in a weak sense via the Wigner transform.
\begin{theorem}[Semi-classical limit of quantum scattering states]\label{thm: semiclassical limit}
Under the assumptions and notations in Theorem \ref{thm: uniform bounds for nonlinear solutions} and Corollary \ref{cor: uniform small-data scattering}, let $f_0\geq0$ be the weak limit of $f_0^{\hbar_j}:= \textup{\textbf{Wig}}^{\hbar_j}[\ga_0^{\hb_j}]$ in $L^2_{q,p}(\mathbb{R}^6)$ for some sequence $\{\hb_j\}_{j=1}^\I \subset (0,1]$ with $\hb_j \to 0$.
\begin{enumerate}[$(i)$]
\item (Compactness and a priori bounds) There exists a non-negative distribution function $g(t)\in C(\mathbb{R}; w-B_{\eta_0})$ such that passing to a subsequence,
$$\lim_{j\to\infty}\mathbf{d}_{[-T,T]}^w\Big(g^{\hbar_j}(t), g(t)\Big)=0$$
for all $T>0$, where
$$g^{\hbar_j}(t):=\mathcal{U}(-t)\textup{\textbf{Wig}}^{\hbar_j}[\ga^{\hb_j}(t)].$$
Moreover, $f(t)=\mathcal{U}(t)g(t)$ satisfies the bounds 
\begin{equation*}
\begin{aligned}
&\|f(t)\|_{L^\I_t(\mathbb{R}; L_{q,p}^1(\mathbb{R}^6))}
+\|\langle p\rangle^{1+a+2\epsilon} f(t)\|_{L^\I_t(\mathbb{R};L_{q,p}^{\frac{3}{2-a-\epsilon}}(\mathbb{R}^6))}\\
&\qquad \qquad \qquad +\|\langle q-tp\rangle^{1+a+2\epsilon} f(t)\|_{L^\I_t(\mathbb{R};L_{q,p}^{\frac{3}{2-a-\epsilon}}(\mathbb{R}^6))}\lesssim \eta_0
\end{aligned}
\end{equation*}
and
$$\sup_{t\in\mathbb{R}}\bigg(\langle t\rangle^{1+a+\epsilon}\|\rho_{f(t)}\|_{L^{\frac{3}{2-a-\epsilon}}(\mathbb{R}^3)}\bigg)\lesssim \eta_0.$$
\item (Semi-classical limit to a global scattering solution) $f(t)$ is a global scattering solution in $C_t(\mathbb{R};H^{-1}_{q,p}(\mathbb{R}^6))$ to the Vlasov equation \eqref{eq: Vlasov} with initial data $f_0$ in the sense of Definition \ref{definition: scattering solutions to Vlasov} such that
$$\lim_{t\to\pm\infty}\|\mathcal{U}(-t)f(t)-f_\pm\|_{H^{-1}_{q,p}(\R^6)}=0.$$
\item (Semi-classical limit of scattering states)
For the scattering states $f_\pm\in H_{q,p}^{-1}(\mathbb{R}^6)$ in $(ii)$, we have
$$\textup{\textbf{Wig}}^{\hbar_j}[\ga^{\hb_j}_\pm]\rightharpoonup f_\pm\mbox{ in }L^2_{q,p}(\mathbb{R}^6).$$
\end{enumerate}
\end{theorem}

\begin{remark}
It is proved in Lions and Paul \cite{Lions Paul 1993} that for a large class of interactions, the Wigner transforms of solutions to the NLH \eqref{eq: NLH} converge weakly-* to a solution to the Vlasov equation \eqref{eq: Vlasov} on a compact interval. 
To the best of the authors' knowledge, Theorem \ref{thm: semiclassical limit} is the first result providing a rigorous justification of the semi-classical limit for \textit{long-time} dynamics. In this direction, many interesting questions remain open, for example, upgrading to the strong convergence as in the recent work of Lafleche and Saffirio \cite{Lafleche Saffirio 2023} as well as including long-range interactions.
\end{remark}

\begin{remark}[Construction of a classical initial distribution $f_0$]\label{rmk:construction of f_0}
By the assumptions of Theorem \ref{thm: uniform bounds for nonlinear solutions} and Lemma \ref{lemma: basic properties of Wigner transform}, $f_0^\hb$ is uniformly bounded in $L^2_{q,p}(\mathbb{R}^6)$. Hence, by the weak compactness of $L^2_{q,p}(\mathbb{R}^6)$ and the Banach-Alaoglu theorem, there exist $\{\hbar_j\}_{j=1}^\I$, with $\hbar_j \to 0$, and non-negative $f_0\in L^2_{q,p}(\mathbb{R}^6)$ such that $f_0^{\hbar_j}$ converges weakly to $f_0$ in $L^2_{q,p}(\mathbb{R}^6)$ and $\|f_0\|_{L_{q,p}^2(\mathbb{R}^6)}\leq \eta_0$.
\end{remark}

\begin{remark}\label{remark: commuting diagram}
The following commutative diagram can express our result.
\[
\begin{tikzcd}[row sep=large, column sep=large]
	\ga_0^\hb \arrow[r, "\textup{solve}"] \arrow[d,->, "\Wig"'] & \ga^\hb(t)\arrow[r, "\textup{scatter}"] \arrow[d,->,"\Wig"'] & \ga_+^\hb \arrow[d,->,"\Wig"'] \\
	f_0^\hb \arrow[d,"\hb\to0"'] & f^\hb(t) \arrow[d,"\hb\to0"'] & f^\hb_+ \arrow[d,"\hb\to0"'] \\
	f_0 \arrow[r,"\textup{solve}"] & f(t) \arrow[r, "\textup{scatter}"] & f_+
	\arrow[from=1-1, to=3-2, phantom, "\circlearrowleft"]
	\arrow[from=1-2, to=3-3, phantom, "\circlearrowleft"]
\end{tikzcd}
\]
Starting from the quantum initial data $\gamma_0^\hbar$, the quantum state $\gamma^\hbar(t)$ evolves via the Hartree equation and subsequently scatters to the quantum states $\gamma_\pm^\hbar$ (from the top left to the top right corner). In parallel, a classical initial datum $f_0$ is obtained as the weak limit of the Wigner transforms of the quantum initial data (from the top left to the bottom left corner). The classical distribution $f(t)$ then evolves via the Vlasov equation and scatters to the classical states $f_\pm$ (from the bottom left to the bottom right corner). Our main result establishes the weak subsequential convergence of the scattering states via the Wigner transform, namely, \( f^{\hbar}_{\pm} = \textup{\textbf{Wig}}^{\hbar}[\gamma^{\hbar}_\pm] \rightharpoonup f_\pm \) (from the top right to the bottom right corner), thereby extending the well-known semiclassical limit result \( f^\hbar(t) = \textup{\textbf{Wig}}^{\hbar}[\gamma^{\hbar}(t)] \rightharpoonup f(t) \) on compact time intervals \cite{Lions Paul 1993}.
\end{remark}

\subsection{New scattering result for the Vlasov equation}

Semi-classical analysis has its own interest, but it sometimes also provides a result that is not easy to obtain directly from standard PDE methods (see \cite{Lewin Sabin 2020}, for example). Indeed, for the Vlasov equation with short-range interaction, small-data scattering is proved based on a fixed-point argument in \cite{Choi Ha 2011, Huang Kwon 2024}, and these results require some regularity and weighted norm bounds on initial data. On the other hand, by our second main result (Theorem \ref{thm: semiclassical limit}), small-data scattering can be deduced without assuming regularity. 

\begin{theorem}[Small-data scattering for the Vlasov equation]\label{thm:scattering Vlasov}
Suppose that $w=\pm\frac{1}{|x|^a}$ with $1<a<\frac{5}{3}$. For any sufficiently small $\epsilon>0$, there exists a small $\eta_0=\eta_0(a,\epsilon)>0$ such that if $f_0\geq0$ and
\begin{equation}\label{eq: initial data condition for f0}
\|f_0\|_{L^1_{q,p}(\mathbb{R}^6)}+\|f_0\|_{L^{\frac{3}{2-a-\epsilon}, 1+a+2\epsilon}_{q,p}(\mathbb{R}^6)}\le \eta_0,
\end{equation}
then there exists a global solution $f(t) \in C(\R_t;w-B_{\eta_0})$ to the Vlasov equation \eqref{eq: Vlasov} with initial data $f(0)=f_0$ in the sense of \eqref{eq:meaning of solution} and scattering states $f_\pm\in H_{q,p}^{-1}(\mathbb{R}^6)$, i.e., $\mathcal{U}(-t)f(t) \to f_\pm$ in $H_{q,p}^{-1}(\mathbb{R}^6)$ as $t\to\pm\infty$. Moreover, $f(t)$ satisfies bounds
\begin{equation}\label{eq:Vlasov density decay}
\|f(t)\|_{C_t(\mathbb{R}; L_{q,p}^1(\mathbb{R}^6))}+\sup_{t \in \R} \bigg(\bt^{1+a+\epsilon}\|\rh_{f(t)}\|_{L^{\frac{3}{2-a-\epsilon}}(\R^3)}\bigg)
\ls \eta_0.
\end{equation}
\end{theorem}

For the proof, we employ the commuting diagram for scattering solutions in the semi-classical limit (see Remark \ref{remark: commuting diagram}). Precisely, for an initial data $f_0$ that satisfies the smallness condition \eqref{eq: initial data condition for f0}, we take its Toeplitz quantization (see \ref{eq: Toeplitz quantization}) to obtain quantum initial data $\gamma_0^\hbar$ such that the smallness condition \eqref{eq: initial data assumption} for quantum states is satisfied and $\textup{\textbf{Wig}}^{\hbar}[\ga_0^{\hb}]\rightharpoonup f_0$ in $L^2_{q,p}(\mathbb{R}^6)$. Then, by Corollary \ref{cor: uniform small-data scattering}, the solution to the NLH \eqref{eq: NLH} with initial data $\gamma_0^\hbar$ scatters. Subsequently, by Theorem \ref{thm: semiclassical limit}, it follows that in the semi-classical limit, the quantum scattering solution converges to a scattering solution to the Vlasov equation with initial data $f_0$.

\begin{remark}
In Theorem \ref{thm:scattering Vlasov}, a small-data scattering solution is constructed, but its uniqueness is not known. Indeed, as long as the density $\rho_f(t)$ is bounded in $L^\infty(\mathbb{R}^3)$, uniqueness is expected among such solutions \cite{Loeper 2006}. However, the bounds in \eqref{eq:Vlasov density decay} are not sufficient for that.
\end{remark}

\subsection{Organization of this paper}
The rest of the paper is divided into two parts; Sections 3 and 4 are for the first main result (Theorem \ref{thm: uniform bounds for nonlinear solutions}), and Sections 5 and 6 are for the second one (Theorem \ref{thm: semiclassical limit}). In Section 3, we show the uniform-in-$\hbar$ dispersive estimates for the free Schr\"odinger flow (Proposition \ref{prop:free dispersive}), and extend it to perturbed linear flows (Proposition \ref{prop:perturbed dispersive}). Then, in Section 4, we establish the uniform decay estimates for small-data solutions to the NLH \eqref{eq: NLH} (Theorem \ref{thm: uniform bounds for nonlinear solutions}) employing the key linear estimates in the previous section. Next, in Section 5, we provide preliminary fundamental tools for semi-classical analysis including the Wigner/Husimi transforms and Weyl/Toeplitz quantizations. Finally, in Section 6, we prove the semi-classical limit of quantum scattering states (Theorem \ref{thm: semiclassical limit}), and a new small-data scattering result of the Vlasov equation (Theorem  \ref{thm:scattering Vlasov}).

\section{Key uniform decay estimates for linear flows}

In this section, we establish key uniform bounds for the free Schr\"odinger equation (Proposition \ref{prop:free dispersive}), and then extend them to perturbed linear flows (Proposition \ref{prop:perturbed dispersive}), where a class of perturbations is chosen to fit into the quantum nonlinear problem \eqref{eq: NLH}.

\subsection{Uniform bounds for the free flow}\label{sec: uniform dispersion estimate for the free flow}
To begin with, we prove density estimates for the quantum free flow.

\begin{proposition}[Uniform bounds for the free Schr\"odinger flow]\label{prop:free dispersive}
For $1\leq r\leq\infty$ and $\sigma>\frac{3}{2r'}$, there exists $C_{\sigma,r}>0$, independent of $\hbar>0$, such that 
\begin{align}
\big\|\rho^\hbar_{\U \gamma_0 \Uinv }\big\|_{L^r_x(\mathbb{R}^3)}&\leq C_{\si,r}\|\langle \hbar\nabla\rangle^\sigma\gamma_0\langle\hbar\nabla\rangle^\sigma\|_{\mathfrak{L}_\hbar^r},\label{eq: free dispersive1}\\
\big\|\rho^\hbar_{\U \gamma_0 \Uinv }\big\|_{L^r_x(\mathbb{R}^3)}&\leq \frac{C_{\si,r}}{|t|^{\frac{3}{r'}}}\|\langle x\rangle^\sigma \gamma_0 \langle x\rangle^\sigma\|_{\mathfrak{L}_\hbar^r}.\label{eq: free dispersive2}
\end{align}
\end{proposition}

\begin{remark}\label{free uniform dispersion estimates remark}
$(i)$ By \eqref{eq: free dispersive2} (resp., \eqref{eq: free dispersive1}) the inequality for orthogonal functions with a gain of summability is obtained. For instance, taking $r=\infty$, $\hbar=1$ and $\gamma_0=\sum_{j=1}^N|\phi_j\rangle\langle\phi_j|$, we have
$$\bigg\|\sum_{j=1}^N|e^{\frac{it}{2}\Delta}\phi_j|^2\bigg\|_{L^\infty(\mathbb{R}^3)}\leq \frac{C_{\sigma,\I}}{(2\pi|t|)^3}\quad\bigg(\textup{resp., }\bigg\|\sum_{j=1}^N|e^{\frac{it}{2}\Delta}\phi_j|^2\bigg\|_{L^\infty(\mathbb{R}^3)}\leq \frac{C_{\sigma,\I}}{(2\pi)^3}\bigg)$$
for an orthonormal set $\{\phi_j\}_{j=1}^N$ in $L^{2,\sigma}(\mathbb{R}^3)=\{\phi: \langle x\rangle^\sigma \phi\in L^2(\mathbb{R}^3)\}$ (resp., $H^\sigma(\mathbb{R}^3)$) for some $\sigma>\frac{3}{2}$.\\
$(ii)$ The inequalities for the Schr\"odinger flow with optimal summability gain are closely linked to those for the transport equation (see \cite[Section 3.3]{Sabin 2014}) via the semi-classical limit. In this context, in the endpoint case $r=\infty$, the dispersion estimate \eqref{eq: free dispersive2} corresponds to the inequality 
$$\|\rho_{\E f_0}\|_{L^\infty(\mathbb{R}^3)}\leq C_{\sigma}\min\Bigg\{\frac{1}{|t|^3}\|\langle q\rangle^{2\sigma}f_0\|_{L^\infty(\mathbb{R}^6)},\|\langle p\rangle^{2\sigma}f_0\|_{L^\infty(\mathbb{R}^6)}\Bigg\},$$
which is a direct consequence of the well-known dispersion estimate 
\begin{equation}\label{eq: dispersion for free transport}
\|\rho_{\E f_0}\|_{L^\infty(\mathbb{R}^3)}\leq \min\Bigg\{\frac{1}{|t|^3}\|f_0\|_{L_q^1(\mathbb{R}^3; L_p^\infty(\mathbb{R}^3))}, \|f_0\|_{L_p^1 (\mathbb{R}^3; L_q^\infty(\mathbb{R}^3))}\Bigg\}
\end{equation}
for the transport equation \cite{Castella Perthame 1996}.\\
$(iii)$ Proposition \ref{prop:free dispersive} is crucially used in this paper, but its proof is elementary (see Section \ref{sec: uniform dispersion estimate for the free flow}), and it implicitly follows from the argument in \cite[Lemma 4]{Lewin Sabin 2014}.
\end{remark}

\begin{proof}[Proof of Proposition \ref{prop:free dispersive}]
Interpolating with the trivial bound
$$\|\rho^\hbar_{\U \gamma_0 \Uinv}\|_{L^1_x}\le(2\pi\hbar)^3\|e^{\frac{it\hbar}{2}\Delta}\gamma_0 e^{-\frac{it\hbar}{2}\Delta}\|_{\mathfrak{S}^1}=(2\pi\hbar)^3\|\gamma_0\|_{\mathfrak{S}^1}=\|\gamma_0\|_{\mathfrak{L}_\hbar^1},$$
it suffices to show the endpoint inequalities:
\begin{align}
\|\rho^\hbar_{\U \gamma_0 \Uinv}\|_{L^\infty_x}
&\leq C_{\sigma, \I} \|\langle \hbar\nabla\rangle^\sigma\gamma_0\langle \hbar\nabla\rangle^\sigma\|_{\mathcal{B}(L^2)},\label{eq: free dispersive1'}\\
\|\rho^\hbar_{\U \gamma_0 \Uinv}\|_{L^\infty_x}
&\leq\frac{C_{\sigma,\I}}{|t|^{3}}\|\langle x\rangle^\sigma \gamma_0 \langle x\rangle^\sigma\|_{\mathcal{B}(L^2)}\label{eq: free dispersive2'}
\end{align}
for $\sigma>\frac{3}{2}$. We will show \eqref{eq: free dispersive1'} and \eqref{eq: free dispersive2'} by duality.

Fix $V\in C_c^\infty(\mathbb{R}^3)$. Then, using the fact that 
$$\int_{\mathbb{R}^3}V(x)\rho^\hb_\gamma(x)dx =(2\pi\hb)^3 \textup{Tr}(V\gamma),$$
where $V$ on the right hand side denotes a multiplication operator, as well as the cyclicity of the trace, we write 
$$\begin{aligned}
\int_{\mathbb{R}^3}V(x)\rho^\hbar_{\U \gamma_0 \Uinv} (x)dx
&=(2\pi\hbar)^3\textup{Tr}\big(V \U \gamma_0 \Uinv \big)\\
&=(2\pi\hbar)^3\textup{Tr}\big(\langle x\rangle^{-\sigma}\Uinv V \U \langle x\rangle^{-\sigma}\langle x\rangle^\sigma\gamma_0\langle x\rangle^\sigma \big).
\end{aligned}$$
Hence, by H\"older's inequality for Schatten class operators, it follows that 
$$\int_{\mathbb{R}^3}V(x)\rho^\hbar_{\U \gamma_0 \Uinv}(x)dx\leq \|\langle x\rangle^{-\sigma}\Uinv \sqrt{|V|}\|_{\mathfrak{L}_\hbar^2}^2\|\langle x\rangle^\sigma \gamma_0 \langle x\rangle^\sigma\|_{\mathcal{B}(L^2)}.$$
Moreover, replacing $\langle x\rangle^\sigma$ by $\langle \hbar\nabla\rangle^\sigma$, one can also show that 
$$\int_{\mathbb{R}^3}V(x)\rho^\hbar_{\U \gamma_0 \Uinv}(x)dx\leq \|\langle \hbar\nabla\rangle^{-\sigma}\Uinv\sqrt{|V|}\|_{\mathfrak{L}_\hbar^2}^2\|\langle \hbar\nabla\rangle^\sigma \gamma_0 \langle \hbar\nabla\rangle^\sigma\|_{\mathcal{B}(L^2)}.$$
Recall that the Schatten $2$-norm $\|\cdot\|_{\mathfrak{S}^2}$ equals to the Hilbert-Schmidt norm, and that $\|T_K\|_{\mathfrak{L}_\hbar^2}=(2\pi\hbar)^{\frac{3}{2}}\|K(x,x')\|_{L_{x,x'}^2}$for an integral operator $T_K$ with kernel $K(x,x')$. Thus, we have
$$\begin{aligned}
\|\langle x\rangle^{-\sigma}\Uinv\sqrt{|V|}\|_{\mathfrak{L}_\hbar^2}^2&=\iint_{\mathbb{R}^{6}} \frac{1}{\langle x\rangle^{2\sigma}}\frac{1}{|t|^3}|V(x')|dxdx'=\frac{C_{\sigma,\infty}}{|t|^3}\|V\|_{L^1},
\end{aligned}$$
where $C_{\sigma,\infty}=\|\langle x\rangle^{-2\sigma}\|_{L_x^1}$, since the kernel of $\Uinv = e^{-\frac{it\hbar}{2}\Delta}$ is given by $\frac{1}{(-2\pi i\hbar t)^{3/2}}e^{\frac{|x-x'|^2}{2i\hbar t}}$. On the other hand, by the Plancherel theorem, we have
$$\begin{aligned}
\|\langle \hbar\nabla\rangle^{-\sigma}\sqrt{|V|}\|_{\mathfrak{L}_\hbar^2}^2&=(2\pi\hbar)^3\|\langle \hbar\nabla\rangle^{-\sigma}\sqrt{|V|}\|_{\mathfrak{S}^2}^2\\
&=(2\pi\hbar)^3\|(\langle \hbar \xi\rangle^{-\sigma})^\vee\|_{L^2}^2\|\sqrt{|V|}\|_{L^2}^2=C_{\sigma,\infty} \|V\|_{L^1}.
\end{aligned}$$
Therefore, \eqref{eq: free dispersive1'} and \eqref{eq: free dispersive2'} are obtained by duality.
\end{proof}

\subsection{Linear Schr\"odinger flows and wave operators}

From Section \ref{sec: notations}, we recall the definition of the vector field $\J_t^\hbar:=x+it\hbar\nabla$ and that of the linear Schr\"odinger flow $\mathcal{U}_V^\hbar(t,t_0)$ associated with the equation $i\hbar\partial_t u=(-\frac{\hbar^2}{2}\Delta +V(t))u$. Then, we define the (finite-time) wave operator
$$\mathcal{W}_V^\hbar(t,t_0):=\Uinv\mathcal{U}_V^\hbar(t,t_0) \CU^\hb(t_0).$$

Note that by the Duhamel formula, the propagator $\mathcal{U}_V^\hbar(t,t_0)$ can be represented as 
\begin{equation}\label{eq: linear Schrodinger, Duhamel}
\mathcal{U}_V^\hbar(t,t_0)
=\CU^\hbar(t-t_0) - \frac{i}{\hb} \int_{t_0}^t \CU^\hb(t-t_1)V(t_1)\mathcal{U}_V^\hbar(t_1,t_0)dt_1.
\end{equation}
By \eqref{eq: linear Schrodinger, Duhamel}, the linear flow $\mathcal{U}_V^\hbar(t,t_0)$ by itself cannot be considered as a perturbation of the free flow $e^{\frac{i(t-t_0)\hbar}{2}\Delta}$ due to the large $\frac{1}{\hbar}$-factor in the integral term. However, in the Heisenberg picture, the equation   $i\hbar\partial_t \gamma=[-\frac{\hbar^2}{2}\Delta+V,\gamma]$ can be considered as a perturbation of the free equation by capturing the hidden $\hbar$-factor in the commutator $[V,\gamma]$. In this sense, the following lemmas will be crucial in our analysis.

\begin{lemma}[Commutator properties of vector fields]\
\begin{enumerate}[$(i)$]
		\item (Commutator with the Schr\"odinger flow)
		\begin{align}
			(i\hbar\nabla)\mathcal{U}_V^\hbar(t,0)&=\mathcal{U}_V^\hbar(t,0)(i\hbar\nabla)+\int_{0}^t \mathcal{U}_V^\hbar(t,t_1)(\nabla V)(t_1)\mathcal{U}_V^\hbar(t_1,0) dt_1,\label{eq: commutator with S flow1}\\
			\J_t^\hbar\mathcal{U}_V^\hbar(t,0)&=\mathcal{U}_V^\hbar(t,0)x+\int_{0}^t \mathcal{U}_V^\hbar(t,t_1)t_1(\nabla V)(t_1)\mathcal{U}_V^\hbar(t_1,0) dt_1.\label{eq: commutator with S flow2}
		\end{align}
		\item (Commutator with the wave operator)
        \begin{align}
			\big[i\hbar\nabla, \mathcal{W}_V^\hbar(t,0)\big]&=\int_0^t \mathcal{W}_V^\hbar(t,t_1)\CU^\hb(t_1)^*(\nabla V)(t_1)\CU^\hb(t_1)\mathcal{W}_V^\hbar(t_1,0) dt_1,\label{eq: commutator with wave operator1}\\
			\big[x, \mathcal{W}_V^\hbar(t,0)\big]&=\int_0^t \mathcal{W}_V^\hbar(t,t_1)\CU^\hb(t_1)^* t_1(\nabla V)(t_1)\CU^\hb(t_1)\mathcal{W}_V^\hbar(t_1,0) dt_1.\label{eq: commutator with wave operator2}
		\end{align}
	\end{enumerate}
\end{lemma}

\begin{proof}
	It suffices to show \eqref{eq: commutator with S flow1} and \eqref{eq: commutator with S flow2}, because they imply \eqref{eq: commutator with wave operator1} and \eqref{eq: commutator with wave operator2} by definitions and \eqref{J to derivative}. Indeed, by direct calculations, one can see that $u(t)=\mathcal{U}_V^\hbar(t,0)\varphi$ satisfies 
	$$i\partial_t(i\hbar\nabla u)=i\hbar\nabla i\partial_t u=i\hbar\nabla\bigg(-\frac{\hbar}{2}\Delta+\frac{1}{\hbar}V\bigg)u=\bigg(-\frac{\hbar}{2}\Delta+\frac{1}{\hbar}V\bigg)(i\hbar\nabla u)+i(\nabla V)u$$
	and
	$$\begin{aligned}
		i\partial_t \big(\J_t^\hbar u\big)&=(x+it\hbar\nabla)i\partial_t u-\hbar\nabla u=(x+it\hbar\nabla)\bigg(-\frac{\hbar}{2}\Delta+\frac{1}{\hbar}V\bigg)u-\hbar\nabla u\\
		&=\bigg(-\frac{\hbar}{2}\Delta+\frac{1}{\hbar}V\bigg)\J_t^\hbar u+it(\nabla V)u.
	\end{aligned}$$
	Hence, by the Duhamel formula, \eqref{eq: commutator with S flow1} and \eqref{eq: commutator with S flow2} follow.
\end{proof}

\subsection{Uniform decay estimate for perturbed flows}
Given $0\le\sigma\leq 2$ and a time interval $I\subset\mathbb{R}$, we define the function space $\mathcal{Z}^\sigma(I)$ as the collection of all time-dependent potentials such that $\|V\|_{\CZ^\si(I)} < \I$, where
\begin{align}\label{eq: Z norm definition}
\|V\|_{\mathcal{Z}^\sigma(I)}:=
   \begin{dcases}
	\|\bt V(t)\|_{L_t^1(I; \dot{W}_x^{1,\infty}(\mathbb{R}^3))} &\text{if } 0\le \si \le 1\\
	\|\langle t\rangle V(t)\|_{L_t^1(I; \dot{W}_x^{1,\infty}(\mathbb{R}^3))}+\|\langle t\rangle V(t)\|_{L_t^1(I; \dot{W}_x^{\sigma,\frac{3}{\sigma-1}}(\mathbb{R}^3))} &\text{if } 1<\si\le 2.
	\end{dcases}
\end{align}
The goal of this subsection is to extend Proposition \ref{prop:free dispersive} to Schr\"odinger flows perturbed by potentials in the class $\mathcal{Z}^\sigma(I)$, which will be a key estimate in Section \ref{sec: uniform bounds for NLH} for a uniform decay estimate for nonlinear solutions.

\begin{proposition}[Uniform dispersion estimate for perturbed flows]\label{prop:perturbed dispersive}
For $3\leq r\leq \infty$ and $\frac{3}{2r'}<\sigma<2$, there exists $K_0=K_0(r,\sigma)\geq1$, independent of $\hbar\in(0,1]$, such that if $\|V\|_{\mathcal{Z}^\sigma(I)}\leq1$ with an interval $I\ni 0$, then for all $t\in I$, 
\begin{equation}\label{eq: perturbed dispersive, small potential}
\|\rho^\hbar_{\mathcal{U}_V^\hbar(t,0)\gamma_0^\hbar \mathcal{U}_V^\hbar(t,0)^*}\|_{L^r_x(\mathbb{R}^3)}\leq\frac{K_0}{\langle t\rangle^{3/r'}}\|\gamma_0^\hbar\|_{\mathfrak{L}^{r,\sigma}}.
\end{equation}
\end{proposition}

For the proof, we make use of the boundedness of the wave operator.

\begin{lemma}[Boundedness of the wave operator]\label{lemma: boundedness of wave operators}
For $0\leq \alpha\leq 2$ and $\al\le \sigma \le 2$, there exists $K_1=K_1(\alpha,\sigma)\geq1$, independent of $\hbar\in(0,1]$, such that if $\|V\|_{\mathcal{Z}^\sigma(I)}\leq1$, then
\begin{align}
\|\langle \hbar\nabla\rangle^\alpha\mathcal{W}_{V}^\hbar(t,0)\langle \hbar\nabla\rangle^{-\alpha}\|_{C_t(I; \mathcal{B}(L^2(\mathbb{R}^3)))}&\leq K_1, \label{eq: wave operator boundedness 1}\\
\|\langle x\rangle^\alpha\mathcal{W}_{V}^\hbar(t,0)\langle x\rangle^{-\alpha}\|_{C_t(I; \mathcal{B}(L^2(\mathbb{R}^3)))}&\leq K_1.\label{eq: wave operator boundedness 2}
\end{align}
\end{lemma}

\begin{proof}
We only prove \eqref{eq: wave operator boundedness 2}, since \eqref{eq: wave operator boundedness 1} can be proved in the same way, except that the commutator relation \eqref{eq: commutator with wave operator1} should be used instead of \eqref{eq: commutator with wave operator2}.

First, we consider the case $\alpha=1$. Note that by the commutator relation \eqref{eq: commutator with wave operator2} for the wave operator and the unitarity of the wave operator,
$$\begin{aligned}
\|[x,\mathcal{W}_V^\hbar(t,0)]\langle x\rangle^{-1}\|_{\mathcal{B}(L^2)}
&\leq \int_I\big\| \mathcal{W}_V^\hbar(t,t')\CU^\hb(t')^* t'(\nabla V)(t') \CU^\hb(t') \mathcal{W}_V^\hbar(t',0) \langle x\rangle^{-1}\big\|_{\mathcal{B}(L^2)} dt'\\
&\leq \int_I|t'|\|V(t')\|_{\dot{W}_x^{1,\infty}} dt'\leq \|V\|_{\mathcal{Z}^\si(I)}\leq 1.
\end{aligned}$$
Hence, together with the unitarity of the wave operator, it follows that 
$$\begin{aligned}
\|\langle x\rangle \mathcal{W}_V^\hbar(t,0)\langle x\rangle^{-1}\|_{\mathcal{B}(L^2)}&\leq\|\mathcal{W}_V^\hbar(t,0)\langle x\rangle^{-1}\|_{\mathcal{B}(L^2)}+\|[x,\mathcal{W}_V^\hbar(t,0)]\langle x\rangle^{-1}\|_{\mathcal{B}(L^2)}\\
&\quad+\|\mathcal{W}_V^\hbar(t,0)x\langle x\rangle^{-1}\|_{\mathcal{B}(L^2)}\lesssim1.
\end{aligned}$$
Then, complex interpolation with the trivial bound $\|\mathcal{W}_V^\hbar(t,0)\|_{\mathcal{B}(L^2)}\leq 1$ yields \eqref{eq: wave operator boundedness 2} with $0\leq\alpha\leq 1$.

Next, we assume that $1\leq\alpha\leq2$. In a similar manner, we write 
$$\begin{aligned}
\|\langle x\rangle^\alpha\mathcal{W}_V^\hbar(t,0)\langle x\rangle^{-\alpha}\|_{\mathcal{B}(L^2)}&\leq\|\langle x\rangle^{\alpha-1}\mathcal{W}_V^\hbar(t,0)\langle x\rangle^{-\alpha}\|_{\mathcal{B}(L^2)}\\
&\quad+\|\langle x\rangle^{\alpha-1}[x,\mathcal{W}_V^\hbar(t,0)]\langle x\rangle^{-\alpha}\|_{\mathcal{B}(L^2)}\\
&\quad+\|\langle x\rangle^{\alpha-1}\mathcal{W}_V^\hbar(t,0)x\langle x\rangle^{-\alpha}\|_{\mathcal{B}(L^2)}\\
&=:\textup{(I)}+\textup{(II)}+\textup{(III)}.
\end{aligned}$$
Then, from \eqref{eq: wave operator boundedness 2} with $0\leq\alpha\leq 1$, it follows that $\textup{(I)}+\textup{(III)}\lesssim1$. For $\textup{(II)}$, applying the commutator relation \eqref{eq: commutator with wave operator2} and \eqref{eq: wave operator boundedness 2} with $0<\alpha-1<1$, we obtain
$$\begin{aligned}
\textup{(II)}
&\leq \int_{-|t|}^{|t|} \big\|\langle x\rangle^{\alpha-1}\mathcal{W}_V^\hbar(t,t') \CU^\hb(t')^* t'(\nabla V)(t') \CU^\hbar(t') \mathcal{W}_V^\hbar(t',0)\langle x\rangle^{-\alpha}\big\|_{\mathcal{B}(L^2)} dt'\\
&\lesssim \int_{-|t|}^{|t|} \big\|\langle x\rangle^{\alpha-1} \CU^\hbar(t')^* t'(\nabla V)(t') \CU^\hbar(t') \langle x\rangle^{-(\alpha-1)}\big\|_{\mathcal{B}(L^2)} dt'\\
&\leq\int_{-|t|}^{|t|} \big\|\langle x\rangle^{\alpha-1} \CU^\hb(t')^* \langle \J_{t'}^\hbar\rangle^{-(\alpha-1)}\big\|_{\mathcal{B}(L^2)}
|t'|\big\|\langle \J_{t'}^\hbar\rangle^{\alpha-1}(\nabla V)(t')\langle \J_{t'}^\hbar\rangle^{-(\alpha-1)}\big\|_{\mathcal{B}(L^2)} \\
&\qquad \quad \cdot 
\big\|\langle \J_{t'}^\hbar\rangle^{\alpha-1} \CU^\hb(t')\langle x\rangle^{-(\alpha-1)}\big\|_{\mathcal{B}(L^2)} dt'.
\end{aligned}$$
On the right hand side, by the representations for the fractional power of $\J_t^\hbar$ in \eqref{eq: factional power of J}, we have
$$\big\|\langle x\rangle^{\alpha-1}\CU^\hbar(t')^* \langle \J_{t'}^\hbar\rangle^{-(\alpha-1)}\big\|_{\mathcal{B}(L^2)}
=\big\|\langle x\rangle^{\alpha-1}\CU^\hbar(t')^* \CU^\hbar(t')\langle x\rangle^{-(\alpha-1)}\CU^\hbar(t')^*\big\|_{\mathcal{B}(L^2)}
=1.$$
Moreover, by \eqref{eq: factional power of J} again, we write 
$$\big\|\langle\J_t^\hbar\rangle^{\alpha-1}\big((\nabla V)(\langle \J_{t}^\hbar\rangle^{-(\alpha-1)}\varphi)\big)\big\|_{L_x^2}
=\Big\|\langle t\hbar\nabla\rangle^{\alpha-1}\Big\{(\nabla V)\Big(\langle t\hbar\nabla\rangle^{-(\alpha-1)}e^{-\frac{i|x|^2}{2t\hbar}}\varphi\Big)\Big\}\Big\|_{L_x^2}.$$
Then, applying the fractional Leibniz rule \cite[Theorem C]{Nahas Ponce 2009} and the Sobolev inequality, we obtain 
$$\begin{aligned}
&\big\|\langle\J_t^\hbar\rangle^{\alpha-1}\big((\nabla V)(\langle \J_{t}^\hbar\rangle^{-(\alpha-1)}\varphi)\big)\big\|_{L_x^2}\\
&\lesssim \|\nabla V\|_{L_x^\infty}\|e^{-\frac{i|x|^2}{2t\hbar}}\varphi\|_{L_x^2}+\big\||t\hbar\nabla|^{\alpha-1}\nabla V\big\|_{L_x^{\frac{3}{\alpha-1}}}\big\|\langle t\hbar\nabla\rangle^{-(\alpha-1)}(e^{-\frac{i|x|^2}{2t\hbar}}\varphi)\big\|_{L_x^{\frac{6}{5-2\alpha}}}\\
&\lesssim \|V\|_{\dot{W}_x^{1,\infty}}\|\varphi\|_{L_x^2}+|t\hbar|^{\alpha-1}\|V\|_{\dot{W}_x^{\alpha, \frac{3}{\alpha-1}}}\big\||\nabla|^{\alpha-1}\langle t\hbar\nabla\rangle^{-(\alpha-1)}(e^{-\frac{i|x|^2}{2t\hbar}}\varphi)\big\|_{L_x^2}\\
&\lesssim \Big(\|V\|_{\dot{W}_x^{1,\infty}}+\|V\|_{ \dot{W}_x^{\alpha,\frac{3}{\alpha-1}}}\Big)\|\varphi\|_{L_x^2},
\end{aligned}$$
where in the last step, we used that $\||\nabla|^{\alpha-1}\langle t\hbar\nabla\rangle^{-(\alpha-1)}\|_{\mathcal{B}(L^2)}=\||\xi|^{\alpha-1}\langle t\hbar\xi\rangle^{-(\alpha-1)}\|_{L_\xi^\infty}\lesssim (|t|\hbar)^{-(\alpha-1)}$. Hence, applying all the above bounds, we obtain
$$\textup{(II)}\lesssim\int_I |t'|\Big(\|V\|_{\dot{W}_x^{1,\infty}}+\|V\|_{ \dot{W}_x^{\alpha,\frac{3}{\alpha-1}}}\Big) dt'\leq \|V\|_{\mathcal{Z}^\sigma(I)}\leq1.$$
Therefore, collecting the estimates for $\textup{(I)}$, $\textup{(II)}$ and $\textup{(III)}$, we prove the desired inequality.
\end{proof}

\begin{proof}[Proof of Proposition \ref{prop:perturbed dispersive}]
Extracting the free propagator, we write
$$\gamma^\hbar(t)=\U \mathcal{W}_V^\hbar(t)\gamma_0^\hbar \mathcal{W}_V^\hbar(t)^*\Uinv,$$
and then applying the dispersion estimate for the free flow (Proposition \ref{prop:free dispersive}), we obtain
$$\|\rho^\hbar_{\gamma^\hbar(t)}\|_{L^r_x}\lesssim \frac{1}{\langle t\rangle^{3/r'}}\|\mathcal{W}_V^\hbar(t)\gamma_0^\hbar \mathcal{W}_V^\hbar(t)^*\|_{\mathfrak{L}_\hbar^{r,\sigma}}.$$
Then, by the boundedness of the wave operator (Lemma \ref{lemma: boundedness of wave operators}),
$$\begin{aligned}
\|\mathcal{W}_V^\hbar(t)\gamma_0^\hbar \mathcal{W}_V^\hbar(t)^*\|_{\mathfrak{L}_\hbar^{r,\sigma}}&\leq \|\langle x\rangle^\sigma \mathcal{W}_V^\hbar(t)\langle x\rangle^{-\sigma}\|_{\mathcal{B}(L^2)}^2
\|\langle x\rangle^\sigma\gamma_0^\hbar \langle x\rangle^\sigma\|_{\mathfrak{L}_\hbar^r} \\
&\quad+\|\sd^\si \mathcal{W}_V^\hbar(t)\sd^{-\sigma}\|_{\mathcal{B}(L^2)}^2
\|\sd^\sigma\gamma_0^\hbar \sd^\sigma\|_{\mathfrak{L}_\hbar^r} \\
&\lesssim\|\gamma_0^\hbar\|_{\mathfrak{L}_\hbar^{r,\sigma}}.
\end{aligned}$$
Therefore, Proposition \ref{prop:perturbed dispersive} follows.
\end{proof}

\section{Uniform bounds for the nonlinear Hartree equation}\label{sec: uniform bounds for the nonlinear Hartree equation}

In this section, we prove our first main result (Theorem \ref{thm: uniform bounds for nonlinear solutions}), that is, the global-in-time uniform bounds for small-data solutions to the nonlinear Hartree equation \eqref{eq: NLH}, employing the uniform dispersion estimates in the previous section. It turns out that our key linear estimate (Proposition \ref{prop:perturbed dispersive}) is somewhat complicated to use directly for the standard contraction mapping argument, in particular, for the difference estimate. To avoid technical difficulty, we divide the proof into two steps. First, in Section \ref{sec: LWP}, we prove a preliminary local well-posedness in a suitable space by the dispersion estimate for the free flow (Proposition \ref{prop:free dispersive}). Second, in Section \ref{sec: uniform bounds for NLH}, we show that the small-data local-in-time solution, constructed in the previous subsection, does not blow up in finite time and it satisfies the desired global-in-time uniform bounds. In addition, in Section \ref{sec: small-data scattering}, we give a proof of small-data scattering (Corollary \ref{cor: uniform small-data scattering}).

\subsection{Local-in-time bound of the density function}\label{sec: LWP}
Throughout this section, we fix $1<a<\frac{5}{3}$ and an arbitrarily small $\epsilon>0$ (depending on $a$), and set 
\begin{equation}\label{eq: choice of r and sigma}
r_\epsilon=\frac{3}{2-a-\epsilon}\quad\textup{and}\quad\sigma_\epsilon=\frac{1+a+2\epsilon}{2}.
\end{equation}
Let $\mathcal{X}_{in}$ denote the collection of all non-negative compact self-adjoint operators such that  
$$\begin{aligned}
\|A\|_{\mathcal{X}_{in}}
:=\|A\|_{\mathfrak{L}_\hbar^1}
+\|\langle \hbar\nabla\rangle^{\sigma_\epsilon}A\langle \hbar\nabla\rangle^{\sigma_\epsilon}\|_{\mathfrak{L}_\hbar^{r_\epsilon}}
\end{aligned}$$
is finite. For solutions, we define $\mathcal{Y}_1(I)$ as the Banach space of time-dependent operators with the norm 
$$\begin{aligned}
\|A\|_{\mathcal{Y}_1(I)}:=\|A\|_{C_t(I;\mathfrak{L}_\hbar^1)}+\|\langle\hbar\nabla\rangle^{\sigma_\epsilon}A\langle\hbar\nabla\rangle^{\sigma_\epsilon}\|_{C_t(I;\mathfrak{L}_\hbar^{r_\epsilon})}
\end{aligned}$$
and $\mathcal{Y}_2(I)$ as that of functions of $(t,x)$ with the norm 
\begin{equation}\label{eq: Y2 norm}
\|\zeta\|_{\mathcal{Y}_2(I)}:=\|\zeta\|_{C_t(I; L_x^1 \cap L^{r_\ep}_x)} .
\end{equation}
Then, $\mathcal{Y}(I)=\mathcal{Y}_1(I)\times\mathcal{Y}_2(I)$ denotes the product space equipped with the norm
$$\|(A,\zeta)\|_{\mathcal{Y}(I)}:=\|A\|_{\mathcal{Y}_1(I)}+\|\zeta\|_{\mathcal{Y}_2(I)}.$$
By the above notations, we formulate the local-in-time bound of the density function for the NLH \eqref{eq: NLH} as follows.

\begin{proposition}[Local well-posedness for the nonlinear Hartree equation]\label{prop: LWP}
Let $1<a<\frac{5}{3}$ and $\epsilon>0$ be sufficiently small.
Then, the following hold:
\begin{enumerate}[$(i)$]
	\item (Existence and uniqueness) For any $\gamma_0^\hbar\in \mathcal{X}_{in}$, there exist $T_\hbar\sim \frac{\hbar}{\|\gamma_0^\hbar\|_{\mathcal{X}_{in}}}$ and a unique solution $\gamma^\hbar(t)\in C_t([-T_\hbar, T_\hbar]; \mathfrak{S}^1)$ to the integral equation 
	\begin{equation}\label{eq: Duhamel form NLH}
		\gamma^\hbar(t)=\U \gamma_0^\hbar \Uinv
        -\frac{i}{\hbar}\int_0^t \CU^\hb(t-t')\big[w*\rho^\hbar_{\gamma^\hbar},\gamma^\hbar\big](t') \CU^\hb(t-t')^* dt'
	\end{equation}
	such that
	\begin{equation}\label{eq: solution bound}
		\|\ga^\hb\|_{\CY_1(I)}\leq 2\|\ga_0^\hb\|_{\CX_{in}}
		\textup{ and }
		\|\rho^\hbar_{\gamma^\hbar}\|_{\mathcal{Y}_2(I)}\leq 2c\|\gamma_0^\hbar\|_{\mathcal{X}_{in}},
	\end{equation}
	where $I=[-T_\hb,T_\hb]$ and $c>0$ is some constant determined independently of $\hbar\in(0,1]$ in the proof.
	\item (Blow-up criteria) Let $(T_{\min}, T_{\max})$ be the maximal existence time interval for the solution $\gamma^\hbar(t)$ constructed in $(i)$. If $T_{\max}>0$ (resp., $T_{\min}<0$) is finite, then   $$\|\langle\hbar\nabla\rangle^{\sigma_\epsilon}\gamma^\hbar(t)\langle\hbar\nabla\rangle^{\sigma_\epsilon}\|_{\mathfrak{L}_\hbar^{r_\epsilon}}\to\infty$$
    as $t\nearrow T_{\max}$ (resp., $t\searrow T_{\min}$).
\end{enumerate}
\end{proposition}

\begin{remark}
In Proposition \ref{prop: LWP} $(i)$, the guaranteed interval of existence depends on $\hbar$, and it may shrink to $\{0\}$ as $\hbar\to 0$. For small-data solutions, it will be shown in the next subsection that their maximal intervals of the solutions are $\mathbb{R}$ and they exist globally in time.
\end{remark}

For the proof, the following lemma is useful.

\begin{lemma}[Interpolation inequalities]\label{lemma: interpolation inequalities}
For $1<a<\frac{5}{3}$ and small $\epsilon>0$, we have
\begin{align}
\|w*\zeta\|_{L^\infty(\mathbb{R}^3)}&\lesssim \|\zeta\|_{L^{\frac{3}{2-a-\epsilon}}(\mathbb{R}^3)}^{\frac{a}{1+a+\epsilon}}\|\zeta\|^{\frac{1+\epsilon}{1+a+\epsilon}}_{L^1(\mathbb{R}^3)},\label{eq: nonlinear estimate 1}\\
\big\||\nabla|^{\frac{1+a+2\epsilon}{2}}(w*\zeta)\big\|_{L^{\frac{6}{1+a+2\epsilon}}(\mathbb{R}^3)}&\lesssim \|\zeta\|_{L^{\frac{3}{2-a-\epsilon}}(\mathbb{R}^3)}^{\frac{a}{1+a+\epsilon}}\|\zeta\|_{L^1(\mathbb{R}^3)}^{\frac{1+\epsilon}{1+a+\epsilon}}.\label{eq: nonlinear estimate 2}
\end{align}
\end{lemma}

\begin{proof}
Decomposing the integral with respect to $R>0$ and applying H\"older's inequality, we obtain
$$|(w*\zeta)(x)|\leq \int_{|x-y|\leq R}+\int_{|x-y|>R}\frac{|\zeta(y)|}{|x-y|^a}dy\lesssim \bigg\|\frac{\mathbbm{1}_{|x|\leq R}}{|x|^a}\bigg\|_{L^{\frac{3}{1+a+\epsilon}}}\|\zeta\|_{L^{\frac{3}{2-a-\epsilon}}}+\frac{1}{R^a}\|\zeta\|_{L^1}.$$
Since $\|\frac{\mathbbm{1}_{|x|\leq R}}{|x|^a}\|_{L^{\frac{3}{1+a+\epsilon}}}\sim R^{1+\epsilon}$, choosing $R>0$ optimizing the upper bound, we prove \eqref{eq: nonlinear estimate 1}. For \eqref{eq: nonlinear estimate 2}, by the Hardy-Littlewood-Sobolev inequality, we obtain 
$$\big\||\nabla|^{\frac{1+a+2\epsilon}{2}} (w*\zeta)\big\|_{L^{\frac{6}{1+a+2\epsilon}}}\lesssim \big\||x|^{-\frac{1+3a+2\epsilon}{2}}*\zeta\big\|_{L^{\frac{6}{1+a+2\epsilon}}}\lesssim\|\zeta\|_{L^{\frac{3}{3-a}}}.
$$
Then, by the interpolation inequality, we prove \eqref{eq: nonlinear estimate 2}.
\end{proof}

Moreover, we need the following lemma to ensure $\rh^\hbar_{\ga^\hb} \in \CY_2(I)$.
\begin{lemma}\label{lem:continuity of density}
	If $\gamma_0^\hbar\in\mathcal{X}_{in}$, then $\rh^\hbar_{\U \ga_0^\hb \Uinv}$ is continuous in $L^r_x(\mathbb{R}^3)$ for any $1\le r \le r_\ep$.
\end{lemma}

\begin{proof}
For $t',t\in\mathbb{R}$, we decompose the difference as 
	\begin{equation*}
		\begin{aligned}
			&\rh^\hbar_{\U \ga_0^\hb \Uinv}(x)- \rh^\hbar_{\CU^\hb(t') \ga_0^\hb \CU^\hb(t')^*}(x) \\
			&=\rh^\hbar_{\CU^\hb(t)(1-\CU^\hb(-t+t'))\ga_0^\hb \CU^\hb(t)^*}(x)
                  +\rh^\hbar_{\CU^\hb(t')\ga_0^\hb (\CU^\hb(t'-t)- 1) \CU^\hb(t')^*}(x)\\
            &=:\textup{(I)}+\textup{(II)}.
		\end{aligned}
	\end{equation*}
For the first term, we observe that by the spectral representation for $\gamma_0^\hbar=\sum_{j=1}^\infty\lambda_j^\hbar|\phi_j^\hbar\rangle\langle\phi_j^\hbar|$ with $\lambda_j^\hbar\geq0$,
$$\begin{aligned}
|\textup{(I)}|^2
&=\Bigg\{\sum_{j=1}^\infty\lambda_j^\hbar\Big(\CU^\hb(t)(1-\CU^\hb(-t+t'))\phi_j^\hbar\Big)(x)\big(\CU^\hb(t)\phi_j^\hbar\big)(x)\Bigg\}^2\\
&\leq \Bigg\{\sum_{j=1}^\infty\lambda_j^\hbar\Big|\Big(\CU^\hb(t)(1-\CU^\hb(-t+t'))\phi_j^\hbar\Big)(x)\Big|^2\Bigg\}\Bigg\{\sum_{j=1}^\infty\lambda_j^\hbar\big|\big(\CU^\hb(t)\phi_j^\hbar\big)(x)\big|^2\Bigg\}\\
&= \rh^\hb_{\CU^\hb(t)(1-\CU^\hb(-t+t')) \ga_0^\hb (1-\CU^\hb(-t+t'))^* \CU^\hb(t)^*}(x)
\rh^\hbar_{\CU^\hb(t) \ga_0^\hb \CU^\hb(t)^* }(x).
\end{aligned}$$
Hence, it follows from \eqref{eq: free dispersive1} that 
\begin{equation*}
\begin{aligned}
\|\textup{(I)}\|_{L^{r}_x}&\le\No{\rh^\hb_{\CU^\hb(t)(1-\CU^\hb(-t+t')) \ga_0^\hb (1-\CU^\hb(-t+t'))^* \CU^\hb(t)^*}}^\tw_{L^r_x}\No{\rh^\hbar_{\CU^\hb(t) \ga_0^\hb  \CU^\hb(t)^* } }_{L^{r}_x}^{\frac{1}{2}}\\
&\le \big \| \sd^{\si} (1-\CU^\hb(-t+t'))\ga_0^\hb (1-\CU^\hb(-t+t'))^*\sd^{\si} \big\|_{\mathfrak{L}_\hbar^{r}}^{\frac{1}{2}}\big \| \sd^{\si}\ga_0^\hb \sd^{\si} \big\|_{\mathfrak{L}_\hbar^{r}}^{\frac{1}{2}},
\end{aligned}
\end{equation*}
provided that $\sigma>\frac{3}{2r'}$ when $r>1$;  $\sigma=0$ when $r=1$. Note here that for $\gamma_0^\hbar\in\mathcal{X}_{in}$, $\| \sd^{\si}\ga_0^\hb \sd^{\si}\|_{\mathfrak{L}_\hbar^{r}}$ is bounded by complex interpolation. Moreover, since $1-\CU^\hb(t-t')^* \to 0$ strongly in $L_x^2$ as $t'\to t$, using the fact $A_n B \to A B$ in $\FS^\al$ if $\al \in [1,\I)$, $A_n \to A$ strongly in $L^2_x$ and $B \in \FS^\al$, we prove that $\|\textup{(I)}\|_{L_x^r}\to 0$ as $t'\to t$. In the same way, one can show that $\|\textup{(II)}\|_{L_x^r}\to 0$. 
\end{proof}

\begin{proof}[Proof of Proposition \ref{prop: LWP}]
Let $I=[-T,T]$ be a short time interval that will be chosen later. By the Banach fixed point theorem, it suffices to show that the nonlinear map $\Gamma^\hbar=(\Gamma_1^\hbar, \Gamma_2^\hbar)$,  given by
$$\left\{\begin{aligned}
\Gamma_1^\hbar(A,\zeta)(t)&:=\U \gamma_0^\hbar \Uinv-\frac{i}{\hbar}\int_0^t \CU^\hb(t-t') [w*\zeta,A](t') \CU^\hb(t-t')^* dt',\\
\Gamma_2^\hbar(A,\zeta)(t)&:=\rho^\hbar_{\Gamma_1^\hbar(A,\zeta)}(t),
\end{aligned}\right.$$
is a contractive map in a ball in the solution space $\mathcal{Y}(I)$. For notational convenience, in the sequel, we omit the time interval $I$ in the norms if there is no confusion.

For $\Gamma_1^\hbar(A,\zeta)$, by the unitarity of the free propagator and the trivial inequality $\|[A,B]\|_{\mathcal{Y}_1}\leq 2\|AB\|_{\mathcal{Y}_1}$, it follows that 
$$\begin{aligned}
\|\Gamma_1^\hbar(A,\zeta)\|_{\mathcal{Y}_1}&\leq \|\gamma_0^\hbar\|_{\mathcal{X}_{in}}+\frac{4T}{\hbar}\bigg\{\|(w*\zeta)A\|_{C_t\mathfrak{L}_\hbar^1}+\|\langle\hbar\nabla\rangle^{\sigma_\epsilon}(w*\zeta)A\langle \hbar\nabla\rangle^{\sigma_\epsilon}\|_{C_t\mathfrak{L}_\hbar^{r_\epsilon}}\bigg\}.
\end{aligned}$$
On the right hand side, by the nonlinear estimate \eqref{eq: nonlinear estimate 1} with $r_\epsilon=\frac{3}{2-a-\epsilon}$, we obtain
\begin{equation}\label{eq: LWP proof 1}
\begin{aligned}
\|(w*\zeta)A\|_{C_t\mathfrak{L}_\hbar^1}&\lesssim \|w*\zeta\|_{L_{t,x}^\infty}\|A\|_{C_t\mathfrak{L}_\hbar^1}\\
& \lesssim  \|\zeta\|_{C_t L_x^{r_\epsilon}}^{\frac{a}{1+a+\epsilon}}
\|\zeta\|_{C_tL_x^1}^{\frac{1+\epsilon}{1+a+\epsilon}}
\|A\|_{C_t\mathfrak{L}_\hbar^1}
\le \|\zeta\|_{\mathcal{Y}_2}\|A\|_{C_t\mathfrak{L}_\hbar^1}.
\end{aligned}
\end{equation}
Moreover, we have
\begin{equation}\label{eq: LWP proof 2}
\begin{aligned}
&\|\langle\hbar\nabla\rangle^{\sigma_\epsilon}(w*\zeta)A \langle\hbar\nabla\rangle^{\sigma_\epsilon}\|_{C_t\mathfrak{L}_\hbar^{r_\epsilon}} \\
&\lesssim  \|\langle\hbar\nabla\rangle^{\sigma_\epsilon}(w*\zeta)\langle\hbar\nabla\rangle^{-\sigma_\epsilon}\|_{C_t\mathcal{B}}\|\langle\hbar\nabla\rangle^{\sigma_\epsilon} A\langle\hbar\nabla\rangle^{\sigma_\epsilon}\|_{C_t\mathfrak{L}_\hbar^{r_\epsilon}}\\
&\lesssim \|\zeta\|_{\mathcal{Y}_2}\|\langle\hbar\nabla\rangle^{\sigma_\epsilon} A\langle\hbar\nabla\rangle^{\sigma_\epsilon}\|_{C_t\mathfrak{L}_\hbar^{r_\epsilon}},
\end{aligned}
\end{equation}
because by the fractional Leibniz rule, the Sobolev inequality and the nonlinear estimates \eqref{eq: nonlinear estimate 1} and \eqref{eq: nonlinear estimate 2}, 
\begin{equation}\label{eq: fractional bilinear estimate}
\begin{aligned}
&\big\|\langle\hbar\nabla\rangle^{\sigma_\epsilon}\big((w*\zeta)\langle\hbar\nabla\rangle^{-\sigma_\epsilon}\phi\big)\big\|_{L^2}\\
&\lesssim \|w*\zeta\|_{L^\infty}\|\phi\|_{L^2}+\hbar^{\frac{1+a+2\epsilon}{2}}\||\nabla|^{\frac{1+a+2\epsilon}{2}} w*\zeta\|_{L^{\frac{6}{1+a+2\epsilon}}}\|\langle\hbar\nabla\rangle^{-\frac{1+a+2\epsilon}{2}}\phi\|_{L^{\frac{6}{2-a-2\epsilon}}}\\
&\lesssim \|\zeta\|_{L^{\frac{3}{2-a-\epsilon}}}^{\frac{a}{1+a+\epsilon}}\|\zeta\|_{L^1}^{\frac{1+\epsilon}{1+a+\epsilon}}\|\phi\|_{L^2}=\|\zeta\|_{L^{r_\epsilon}}^{\frac{a}{1+a+\epsilon}}\|\zeta\|_{L^1}^{\frac{1+\epsilon}{1+a+\epsilon}}\|\phi\|_{L^2}.
\end{aligned}
\end{equation}
Therefore, we prove that 
\begin{equation}\label{eq: Gamma1 bound}
\|\Gamma_1^\hbar(A,\zeta)\|_{\mathcal{Y}_1}\leq \|\gamma_0^\hbar\|_{\mathcal{X}_{in}}+\frac{cT}{\hbar}\|\zeta\|_{\mathcal{Y}_2}\|A\|_{\mathcal{Y}_1}.
\end{equation}
For $\Gamma_2^\hbar(A,\zeta)$, by the fact $\|\rho_A\|_{L_x^1}\leq\|A\|_{\mathfrak{S}^1}$ and applying the  bound for the free evolution (Proposition \ref{prop:free dispersive}) with $\sigma_\epsilon=\frac{1+a+2\epsilon}{2}>\frac{3}{2r_\epsilon'}=\frac{1+a+\epsilon}{2}$ and small $T>0$, we obtain
$$\begin{aligned}
\|\Gamma_2^\hbar(A,\zeta)\|_{\mathcal{Y}_2}&\leq \|\Gamma_1^\hbar(A,\zeta)\|_{C_t\mathfrak{L}_\hbar^1}+\sup_{t\in I}\|\Gamma_2^\hbar(A,\zeta)(t)\|_{L^{r_\epsilon}}\\
&\lesssim \|\gamma_0^\hbar\|_{\mathcal{X}_{in}}+\frac{2}{\hbar}\Big(\|(w*\zeta)A\|_{L_t^1\mathfrak{L}_\hbar^1}+\|\langle\hbar\nabla\rangle^{\sigma_\epsilon}(w*\zeta)A\langle \hbar\nabla\rangle^{\sigma_\epsilon}\|_{L_t^1\mathfrak{L}^{r_\epsilon}_\hb}\Big).
\end{aligned}$$
Then, it follows from \eqref{eq: LWP proof 1} and \eqref{eq: LWP proof 2} that 
\begin{equation}\label{eq: Gamma2 bound}
\|\Gamma_2^\hbar(A,\zeta)\|_{\mathcal{Y}_2}\leq c\|\gamma_0^\hbar\|_{\mathcal{X}_{in}}+\frac{cT}{\hbar}\|\zeta\|_{\mathcal{Y}_2}\|A\|_{\mathcal{Y}_1}.
\end{equation}
Now, we take $T=\frac{\hbar}{100\max\{c,1\}^2 \|\gamma_0^\hbar\|_{\mathcal{X}_{in}}}$ and set
$$\wt{\mathcal{Y}}(I):=\Big\{(A,\zeta)\in \mathcal{Y}_1(I)\times \mathcal{Y}_2(I):\ \|A\|_{\mathcal{Y}_1}\leq 2\|\gamma_0^\hbar\|_{\mathcal{X}_{in}}\textup{ and }\|\zeta\|_{\mathcal{Y}_2}\leq 2c\|\gamma_0^\hbar\|_{\mathcal{X}_{in}}\Big\},$$
where $I=[-T,T]$.
Then, it follows from \eqref{eq: Gamma1 bound} and \eqref{eq: Gamma2 bound} that $\Gamma$ maps from $\wt{\mathcal{Y}}(I)$ to itself. 

Repeating the arguments in the proof of \eqref{eq: Gamma1 bound} and \eqref{eq: Gamma2 bound}, one can show that the differences 
$$\left\{\begin{aligned}
\Gamma_1^\hbar(A,\zeta)(t)-\Gamma_1^\hbar(\tilde{A},\tilde{\zeta})(t)
&=-\frac{i}{\hbar}\int_0^t \CU^\hb(t-t') \big([w*\zeta,A]-[w*\tilde{\zeta},\tilde{A}]\big)(t') \CU^\hb(t-t')^* dt',\\
\Gamma_2^\hbar(A,\zeta)(t)-\Gamma_2^\hbar(\tilde{A},\tilde{\zeta})(t)
&=-\frac{i}{\hbar}\int_0^t \rho^\hbar_{\CU^\hb(t-t')([w*\zeta,A]-[w*\tilde{\zeta},\tilde{A}])(t') \CU^\hb(t-t')^*} dt'
\end{aligned}\right.$$
satisfy that if $(A,\zeta), (\tilde{A},\tilde{\zeta})\in \mathcal{Y}$, then $\|\Gamma^\hbar(A,\zeta)-\Gamma^\hbar(\tilde{A},\tilde{\zeta})\|_{\mathcal{Y}}\leq \frac{1}{2}\|(A,\zeta)-(\tilde{A},\tilde{\zeta})\|_{\mathcal{Y}}$. Therefore, we conclude that $\Gamma^\hbar$ is contractive in $\mathcal{Y}$. 

For $(ii)$, by construction, we immediately obtain the blow-up criteria, 
$$\begin{aligned}
\|\gamma^\hbar(t)\|_{\mathfrak{L}_\hbar^1}+\|\langle\hbar\nabla\rangle^{\sigma_\epsilon}\gamma^\hbar(t)\langle\hbar\nabla\rangle^{\sigma_\epsilon}\|_{\mathfrak{L}_\hbar^{r_\epsilon}}\to\infty
\end{aligned}$$
as $t\nearrow T_{\max}$ (resp., $t\searrow T_{\min}$). However, in the above, $\|\gamma^\hbar(t)\|_{\mathfrak{L}_\hbar^1}$ can be removed by the conservation law of the total number of particles, i.e., $\textup{Tr}(\gamma^\hbar(t))=\textup{Tr}(\gamma_0^\hbar)$.
\end{proof}

\subsection{Proof of uniform bounds for the nonlinear Hartree equation (Theorem \ref{thm: uniform bounds for nonlinear solutions})}\label{sec: uniform bounds for NLH}

Next, we prove the uniform decay bound \eqref{eq: dispersion estimate} using the standard bootstrap argument. 

\subsubsection{Bootstrap argument for the density function}
For sufficiently small $\eta_0>0$, we assume that initial data $\gamma_0^\hbar$ satisfies the assumptions in Theorem \ref{thm: uniform bounds for nonlinear solutions}, that is,
$$\|\gamma^\hbar_0\|_{\mathfrak{L}_\hbar^1}+\|\gamma^\hbar_0\|_{\mathfrak{L}_\hbar^{r_\epsilon,\sigma_\epsilon}}\leq\eta_0,$$
where $r_\epsilon=\frac{3}{2-a-\epsilon}$ and $\sigma_\epsilon=\frac{1+a+2\epsilon}{2}$, and let $\gamma^\hbar(t)$
be the solution to the nonlinear Hartree equation \eqref{eq: NLH}, constructed in Proposition \ref{prop: LWP}, such that 
\begin{equation}\label{eq: primitive bound for density on a short time interval}
\|\rho^\hbar(t)\|_{C_t([-T_\hbar, T_\hbar]; L_x^1(\mathbb{R}^3)\cap L^{r_\ep}_x(\mathbb{R}^3))} \le 2c\eta_0,
\end{equation}
where $\rho^\hbar(t):=\rho^\hbar_{\gamma^\hbar(t)}$ and $T_\hbar>0$ is a sufficiently small number which may depend on $\hbar\in(0,1]$. The following lemma shows that the desired potential bounds can be obtained from the density bounds.
\begin{lemma}[Interpolation inequalities for potentials]\label{lemma: interpolation inequalities for potentials}
$$\|\nabla (w*\zeta)\|_{L^\infty(\mathbb{R}^3)}+\|w*\zeta\|_{\dot{W}_x^{\frac{1+a+2\epsilon}{2},\frac{6}{a-1+2\epsilon}}(\mathbb{R}^3)}\lesssim \|\zeta\|_{L^{\frac{3}{2-a-\epsilon}}_\epsilon(\mathbb{R}^3)}^{\frac{1+a}{1+a+\epsilon}}\|\zeta\|_{L^1(\mathbb{R}^3)}^{\frac{\epsilon}{1+a+\epsilon}}.$$
\end{lemma}

\begin{proof}
Repeating the argument to prove the interpolation inequalities (Lemma \ref{lemma: interpolation inequalities}), we prove that 
$$\begin{aligned}
|\nabla(w*\zeta)(x)|&\leq a\int_{|x-y|\leq R}+\int_{|x-y|\geq R}\frac{\zeta(y)}{|x-y|^{a+1}}dy\\
&\lesssim R^{\epsilon}\|\zeta\|_{L^{\frac{3}{2-a-\epsilon}}}+\frac{1}{R^{a+1}}\|\zeta\|_{L^1}\sim \|\zeta\|_{L^{\frac{3}{2-a-\epsilon}}}^{\frac{1+a}{1+a+\epsilon}}\|\zeta\|_{L^1}^{\frac{\epsilon}{1+a+\epsilon}},
\end{aligned}$$
where in the last step, the optimizing $R>0$ is chosen. 
On the other hand, by the Hardy-Littlewood-Sobolev inequality with  $|\nabla|^{\frac{1+a+2\epsilon}{2}}|\nabla|^{-(3-a)}=|\nabla|^{-(3-\frac{1+3a+2\epsilon}{2})}$ and the interpolation inequality, we obtain that 
\begin{equation}\label{eq: why a<5/3}
\|w*\zeta\|_{\dot{W}_x^{\frac{1+a+2\epsilon}{2},\frac{6}{a-1+2\epsilon}}}\lesssim\big\||x|^{-\frac{1+3a+2\epsilon}{2}}*\zeta\big\|_{L_x^{\frac{6}{a-1+2\epsilon}}}\lesssim\|\zeta\|_{L_x^{\frac{3}{2-a}}}\leq \|\zeta\|_{L^1}^{\frac{\epsilon}{1+a+\epsilon}}\|\zeta\|_{L^{\frac{3}{2-a-\epsilon}}}^{\frac{1+a}{1+a+\epsilon}}.
\end{equation}
Note here that in the second inequality, the assumption $a<\frac{5}{3}$ is crucially used in \eqref{eq: why a<5/3}, because we need  $0<\frac{1+3a+2\epsilon}{2}<3$ for the Hardy-Littlewood-Sobolev inequality.
\end{proof}

Recalling the definition of the $\mathcal{Z}^\sigma(I)$-norm (see \eqref{eq: Z norm definition}), we note that by Lemma \ref{lemma: interpolation inequalities for potentials} and \eqref{eq: primitive bound for density on a short time interval} with $\frac{1+a+2\epsilon}{2}=\sigma_\epsilon$, $\frac{6}{a-1+2\epsilon}=\frac{3}{\sigma_\epsilon-1}$ and the choice of sufficiently small $\eta_0>0$, we have
\begin{equation}\label{eq: first potential Z norm bound}
\|\Phi^\hbar\|_{\mathcal{Z}^{\sigma_\epsilon}([-T_\hbar, T_\hbar])}\leq 1.
\end{equation}
Subsequently, the dispersion estimate (Proposition \ref{prop:perturbed dispersive}) yields 
$$\|\rh^\hbar(t)\|_{\wt{\CY}_2([-T_\hbar, T_\hbar])}\leq(K_0+1)\eta_0,$$
where $K_0\geq1$ is the constant given in Proposition \ref{prop:perturbed dispersive} and 
$$\|\rh^\hbar(t)\|_{\wt{\CY}_2(I)}:= \|\rh^\hbar(t)\|_{C_t(I;L^1_x(\mathbb{R}^3))}
+\sup_{t\in I} \bt^{\frac{3}{r_\ep'}}\|\rho^\hbar(t)\|_{L^{r_\ep}_x(\mathbb{R}^3)}.$$
Now, we let $I_*=(T_*^-, T_*^+)$ denote the maximal time interval such that $I_*\ni 0$ and
\begin{equation}\label{eq: a priori bound}
\|\rh^\hbar(t)\|_{\wt{\CY}_2(I_*)}\leq 2(K_0+1)\eta_0\ll1.
\end{equation}
Then, for \eqref{eq: dispersion estimate}, it suffices to show $I_*=\mathbb{R}$. Indeed, repeating the proof of \eqref{eq: first potential Z norm bound}, one can use Lemma \ref{lemma: interpolation inequalities for potentials} to show that $\|\Phi^\hbar(t)\|_{\mathcal{Z}^{\sigma_\epsilon}(I_*)}\leq 1$. Moreover, $\ga^\hb(t)$ obeys\footnote{One can check that the two formulations \eqref{eq: Duhamel form NLH} and \eqref{eq: linear propagator form NLH} are equivalent by simple algebraic calculation.}
\begin{equation}\label{eq: linear propagator form NLH}
\gamma^\hbar(t)=\mathcal{U}_{\Phi^\hbar}^\hbar(t,0)\gamma_0^\hbar\mathcal{U}_{\Phi^\hbar}^\hbar(t,0)^*.
\end{equation}
Hence, it follows from the dispersion estimate (Proposition \ref{prop:perturbed dispersive}) that $\|\rh^\hbar(t)\|_{\wt{\CY}_2(I_*)}\leq(K_0+1)\eta_0$. If $T_*^+$ or $T_*^-$ is finite, then by \eqref{eq: a priori bound}, it contradicts the maximality of the interval $I_*$. We conclude that \eqref{eq: dispersion estimate} holds, but \eqref{eq: dispersion estimate'} is also satisfied by Lemma \ref{lemma: interpolation inequalities for potentials}.

\subsubsection{Estimates for the density operator}
Applying the boundedness of the wave operator (Lemma \ref{lemma: boundedness of wave operators}) to the nonlinear solution 
$$\gamma^\hbar(t)=\mathcal{U}_{\Phi^\hbar}^\hbar(t,0)\gamma_0^\hbar \mathcal{U}_{\Phi^\hbar}^\hbar(t,0)^*=\mathcal{U}^\hbar(t)\mathcal{W}_{\Phi^\hbar}^\hbar(t,0)\gamma_0^\hbar \mathcal{W}_{\Phi^\hbar}^\hbar(t,0)^*\mathcal{U}^\hbar(t)^*,$$
we obtain that for all $t\in\mathbb{R}$,
$$\begin{aligned}
\|\langle\hbar\nabla\rangle^{\sigma_\epsilon}\gamma^\hbar(t)\langle\hbar\nabla\rangle^{\sigma_\epsilon}\|_{\mathfrak{L}_\hbar^{r_\epsilon}}&\leq\|\langle\hbar\nabla\rangle^{\sigma_\epsilon}\mathcal{W}_{\Phi^\hbar}^\hbar(t,0)\langle\hbar\nabla\rangle^{-\sigma_\epsilon}\|_{\mathcal{B}(L^2)}^2\|\langle\hbar\nabla\rangle^{\sigma_\epsilon}\gamma_0^\hbar\langle\hbar\nabla\rangle^{\sigma_\epsilon}\|_{\mathfrak{L}_\hbar^{r_\epsilon}}\\
&\leq K_1^2\|\langle\hbar\nabla\rangle^{\sigma_\epsilon}\gamma_0^\hbar\langle\hbar\nabla\rangle^{\sigma_\epsilon}\|_{\mathfrak{L}^{r_\epsilon}},\\
\|\langle \J_t^\hbar\rangle^{\sigma_\epsilon}\gamma^\hbar(t)\langle \J_t^\hbar\rangle^{\sigma_\epsilon}\|_{\mathfrak{L}_\hbar^{r_\epsilon}}&\leq\|\langle x\rangle^{\sigma_\epsilon}\mathcal{W}_{\Phi^\hbar}^\hbar(t,0)\langle x\rangle^{-\sigma_\epsilon}\|_{\mathcal{B}(L^2)}^2\|\langle x\rangle^{\sigma_\epsilon}\gamma_0^\hbar\langle x\rangle^{\sigma_\epsilon}\|_{\mathfrak{L}_\hbar^{r_\epsilon}}\\
&\leq K_1^2 \|\langle x\rangle^{\sigma_\epsilon}\gamma_0^\hbar\langle x\rangle^{\sigma_\epsilon}\|_{\FL_\hbar^{r_\epsilon}}.
\end{aligned}$$
Therefore, it follows that 
\begin{equation}\label{eq: a posteriori operator bound}
	\begin{aligned}
\big\|\langle\hbar\nabla\rangle^{\sigma_\epsilon}\gamma^\hbar(t)\langle\hbar\nabla\rangle^{\sigma_\epsilon}\big\|_{C_t(\mathbb{R}; \mathfrak{L}_\hbar^{r_\epsilon})}+\big\|\langle \J_t^\hbar\rangle^{\sigma_\epsilon}\gamma^\hbar(t)\langle \J_t^\hbar\rangle^{\sigma_\epsilon}\big\|_{C_t(\mathbb{R}; \mathfrak{L}_\hbar^{r_\epsilon})}
\le (1+2K_1^2)\eta_0.
\end{aligned}
\end{equation}

\subsection{Improved uniform dispersion estimate (Proof of Proposition \ref{prop: upgraded uniform bounds})}

For Proposition \ref{prop: upgraded uniform bounds} $(i)$, we assume that $1<a<\frac{3}{2}$. Then, by Proposition \ref{prop:perturbed dispersive} with $r=\infty$ and $\sigma=\frac{3}{2}+\epsilon$, it suffices to show that $\|\Phi^\hbar\|_{\mathcal{Z}^{\frac{3}{2}+\epsilon}(\mathbb{R})}\leq 1$ (see \eqref{eq: Z norm definition} for the definition of the $\mathcal{Z}^\sigma(I)$-norm). Indeed, by the Sobolev inequality, the interpolation inequality and Theorem \ref{thm: uniform bounds for nonlinear solutions}, we obtain
\begin{equation}\label{eq: gradient potential bound proof'}
\begin{aligned}
\|\Phi^\hbar(t)\|_{\dot{W}^{\frac{3}{2}+\epsilon,\frac{6}{1+2\epsilon}}(\mathbb{R}^3)}&\lesssim \big\||x|^{-(a+\frac{3}{2}+\epsilon)}*\rho^\hbar(t)\big\|_{L^{\frac{6}{1+2\epsilon}}(\mathbb{R}^3)}\lesssim \|\rho^\hbar(t)\|_{L^{\frac{3}{2-a}}(\mathbb{R}^3)}\\
&\leq \|\rho^\hbar(t)\|_{L^{\frac{3}{2-a-\epsilon}}(\mathbb{R}^3)}^{\frac{1+a}{1+a+\epsilon}}\|\rho^\hbar(t)\|_{L^1(\mathbb{R}^3)}^{\frac{\epsilon}{1+a+\epsilon}}\\
&\lesssim\bigg(\frac{\eta_0}{\langle t\rangle^{1+a+\epsilon}}\bigg)^{\frac{1+a}{1+a+\epsilon}}\eta_0^{\frac{\epsilon}{1+a+\epsilon}}=\frac{\eta_0}{\langle t\rangle^{1+a}}.
\end{aligned}
\end{equation}
On the other hand, by Lemma \ref{lemma: interpolation inequalities for potentials}, we have
$$\|\nabla\Phi^\hbar(t)\|_{L^\infty(\mathbb{R}^3)}\lesssim\|\rho^\hbar(t)\|_{L^{\frac{3}{2-a-\epsilon}}(\mathbb{R}^3)}^{\frac{1+a}{1+a+\epsilon}}\|\rho^\hbar(t)\|_{L^1(\mathbb{R}^3)}^{\frac{\epsilon}{1+a+\epsilon}}\lesssim\frac{\eta_0}{\langle t\rangle^{1+a}}.$$
Therefore, it follows that $\|\Phi^\hbar\|_{\mathcal{Z}^{\frac{3}{2}+\epsilon}(\mathbb{R})}\leq 1$.

For Proposition \ref{prop: upgraded uniform bounds} $(ii)$, by Proposition \ref{prop:perturbed dispersive} with $r=\frac{3}{2a-3+2\epsilon}$ and $\sigma=\frac{6-2a-\epsilon}{2}$, it is enough to show that 
$$\|\Phi^\hbar\|_{\mathcal{Z}^{\frac{6-2a-\epsilon}{2}}(\mathbb{R})}\leq 1 $$
(see \eqref{eq: Z norm definition}), provided that $\frac{3}{2}\leq a<\frac{5}{3}$. Indeed, repeating the estimates in \eqref{eq: gradient potential bound proof'}, one can show that 
, because 
$$\begin{aligned}
\|\Phi^\hbar(t)\|_{\dot{W}^{\frac{6-2a-\epsilon}{2},\frac{6}{4-2a-\epsilon}}(\mathbb{R}^3)}&\lesssim \big\||x|^{-\frac{6-\epsilon}{2}}*\rho^\hbar_{\gamma^\hbar(t)}\big\|_{L^{\frac{6}{4-2a-\epsilon}}(\mathbb{R}^3)}\lesssim\|\rho^\hbar(t)\|_{L^{\frac{3}{2-a}}(\mathbb{R}^3)}\\
&\lesssim\|\rho^\hbar(t)\|_{L^{\frac{3}{2-a-\epsilon}}(\mathbb{R}^3)}^{\frac{1+a}{1+a+\epsilon}}\|\rho^\hbar(t)\|_{L^1(\mathbb{R}^3)}^{\frac{\epsilon}{1+a+\epsilon}}\lesssim \frac{\eta_0}{\langle t\rangle^{1+a}}
\end{aligned}$$
and $\|\nabla\Phi^\hbar(t)\|_{L^\infty(\mathbb{R}^3)}\lesssim\frac{\eta_0}{\langle t\rangle^{1+a}}$.
Thus, it follows that $\|\Phi^\hbar\|_{\mathcal{Z}^{\frac{6-2a-\epsilon}{2}}(\mathbb{R})}\leq 1 $.

\subsection{Proof of small-data scattering (Proof of Corollary \ref{cor: uniform small-data scattering})}\label{sec: small-data scattering}

For $t_2>t_1$, by the Duhamel formula \eqref{eq: crude Duhamel formula}, we write 
$$\mathcal{U}^\hbar(t_2)^*\ga^\hb(t_2) \mathcal{U}^\hbar(t_2)-\mathcal{U}^\hbar(t_1)^* \ga^\hb(t_1) \mathcal{U}^\hbar(t_1)
=-\frac{i}{\hbar}\int_{t_1}^{t_2} \mathcal{U}^\hbar(t')^*\big[\Phi^\hb,\ga^\hb\big](t') \mathcal{U}^\hbar(t') dt',$$
where $\mathcal{U}^\hbar(t)=e^{\frac{it\hbar}{2}\Delta}$. Then, it follows from the unitarity of the free propagator and the trivial inequality $\|[A,B]\|_{\mathfrak{L}_\hbar^1}\leq 2\|AB\|_{\mathfrak{L}_\hbar^1}$ for self-adjoint operators $A$ and $B$ that 
$$\big\|\mathcal{U}^\hbar(t_2)^*\ga^\hb(t_2) \mathcal{U}^\hbar(t_2)-\mathcal{U}^\hbar(t_1)^* \ga^\hb(t_1) \mathcal{U}^\hbar(t_1)
\big\|_{\mathfrak{L}_\hbar^1}
\leq\frac{2}{\hb}\int_{t_1}^{t_2}\|\Phi^\hbar(t')\|_{L^\infty(\mathbb{R}^3)}\|\gamma^\hbar(t')\|_{\mathfrak{L}_\hbar^1}dt'.$$
Note that by the conservation law, $$\|\gamma^\hbar(t')\|_{\mathfrak{L}_\hbar^1}=\|\gamma_0^\hbar\|_{\mathfrak{L}_\hbar^1}\leq\eta_0.$$
On the other hand, by the interpolation inequality \eqref{eq: nonlinear estimate 1} and the uniform bounds (Theorem \ref{thm: uniform bounds for nonlinear solutions}), we have
$$\|\Phi^\hbar(t')\|_{L_x^\infty(\mathbb{R}^3)}\lesssim \|\rho^\hbar(t')\|_{L_x^{\frac{3}{2-a-\epsilon}}(\mathbb{R}^3)}^{\frac{a}{1+a+\epsilon}}\|\rho^\hbar(t')\|_{L_x^1(\mathbb{R}^3)}^{\frac{1+\epsilon}{1+a+\epsilon}}\lesssim \bigg(\frac{\eta_0}{\langle t'\rangle^{1+a+\epsilon}}\bigg)^{\frac{a}{1+a+\epsilon}}\eta_0^{\frac{1+\epsilon}{1+a+\epsilon}}=\frac{\eta_0}{\langle t'\rangle^a}.$$
Thus, by the short-range assumption $a>1$, it follows that 
\begin{equation}\begin{aligned}\label{eq: convergence to the scattering state}
\big\|\mathcal{U}^\hbar(t_2)^*\ga^\hb(t_2) \mathcal{U}^\hbar(t_2)-\mathcal{U}^\hbar(t_1)^* \ga^\hb(t_1) \mathcal{U}^\hbar(t_1) \big\|_{\mathfrak{L}_\hbar^1}
&\lesssim\frac{\eta_0^2}{\hb}\int_{t_1}^{t_2}\frac{dt'}{\langle t'\rangle^a} \\
&\ls \frac{\eta_0^2}{\hb \langle \min(|t_1|,|t_2|)\rangle ^{a-1}}\to 0
\end{aligned}\end{equation}
as $t_2 \to -\infty$ (resp., as $t_1\to+\infty$). This shows that the limits $\displaystyle\gamma_\pm^\hbar:=\lim_{t\to\pm\infty}\mathcal{U}^\hbar(t)^*\ga^\hb(t) \mathcal{U}^\hbar(t)$ exist in $\mathfrak{L}_\hbar^1$.

\section{Preliminaries for semi-classical analysis}\label{sec: Preliminaries for semi-classical analysis}

In this section, we summarize the basic properties of the Wigner transform and the Weyl quantization as well as the Husimi transform and the Toeplitz quantization, which will be used in the next section for semi-classical analysis of quantum scattering states. For more details, we refer to Lions-Paul \cite[Section III]{Lions Paul 1993} and Lafleche \cite{Lafleche 2024}. In addition, we derive the Wigner transform of the linear Schr\"odinger equation.

\subsection{Wigner transform, Weyl quantization, and their variants}\label{sec: Wigner transform and Weyl quantization}

The \textit{Weyl quantization} $$\mathbf{Weyl}^\hbar: L_{q,p}^2(\mathbb{R}^{6})\to \mathfrak{L}_\hbar^2$$
is a map which sends functions on the phase space to operators. Precisely, for a distribution function $f\in L_{q,p}^2(\mathbb{R}^{6})$, its Weyl quantization $\mathbf{Weyl}^\hbar[f]$ is defined as the integral operator with kernel 
\begin{equation}\label{eq: Weyl quantization}
\begin{aligned}
\mathbf{Weyl}^\hbar[f](x,x'):&=\frac{1}{(2\pi\hbar)^3}\int_{\mathbb{R}^3}f\bigg(\frac{x+x'}{2},p\bigg) e^{\frac{ip\cdot(x-x')}{\hbar}}dp\\
&=\frac{1}{\hbar^3}(\mathcal{F}_p^{-1}f)\bigg(\frac{x+x'}{2},\frac{x-x'}{\hbar}\bigg).
\end{aligned}\end{equation}
On the other hand, the \textit{Wigner transform} $$\mathbf{Wig}^\hbar: \mathfrak{L}_\hbar^2\to L_{q,p}^2(\mathbb{R}^{6})$$
is the inverse map defined by 
\begin{equation}\label{eq: Wigner transform}
\mathbf{Wig}^\hbar[\gamma](q,p):=\int_{\mathbb{R}^3}\gamma\bigg(q+\frac{y}{2},q-\frac{y}{2}\bigg) e^{-\frac{iy\cdot p}{\hbar}}dy=\mathcal{F}_y\bigg(\gamma\bigg(q+\frac{\cdot}{2},q-\frac{\cdot}{2}\bigg)\bigg)\bigg(\frac{p}{\hbar}\bigg),
\end{equation}
where $\gamma(x,x')$ denotes the integral kernel of $\gamma$. These two transforms are isometric isomorphisms, and they are inverses of each other.
\begin{lemma}[Basic properties of the Wigner transform and the Weyl quantization]\label{lemma: basic properties of Wigner transform}\
\begin{enumerate}[$(i)$]
\item $\|\mathbf{Wig}^\hbar[\gamma]\|_{L_{q,p}^2(\mathbb{R}^{6})}=\|\gamma\|_{\mathfrak{L}_\hbar^2}$ and $\|\mathbf{Weyl}^\hbar[f]\|_{\mathfrak{L}_\hbar^2}=\|f\|_{L_{q,p}^2(\mathbb{R}^{6})}$.
\item $\mathbf{Wig}^\hbar\circ\mathbf{Weyl}^\hbar=\textup{Id}_{L_{q,p}^2(\mathbb{R}^{6})}$ and $\mathbf{Weyl}^\hbar\circ\mathbf{Wig}^\hbar=\textup{Id}_{\mathfrak{L}_\hbar^2}$.
\end{enumerate}
\end{lemma}

\begin{proof}
By \eqref{eq: Weyl quantization} and \eqref{eq: Wigner transform}, the lemma follows from the Plancherel theorem and the fact that $\|\gamma\|_{\mathfrak{L}_\hbar^2}=(2\pi\hbar)^{\frac{3}{2}}\|\gamma(x,x')\|_{L_{x,x'}^2}$.
\end{proof}

\begin{remark}\label{remark: Wigner Weyl remark}
By polarization of Lemma \ref{lemma: basic properties of Wigner transform} $(i)$,
$$\begin{aligned}
\big\langle \mathbf{Wig}^\hbar[\gamma_1]\big|\mathbf{Wig}^\hbar[\gamma_2]\big\rangle_{L_{q,p}^2(\mathbb{R}^6)}&=(2\pi\hbar)^3\textup{Tr}(\gamma_1^*\gamma_2),\\
(2\pi\hbar)^3\textup{Tr}\big(\mathbf{Weyl}^\hbar[f_1]^*\mathbf{Weyl}^\hbar[f_2]\big)&=\langle f_1|f_2\rangle_{L_{q,p}^2(\mathbb{R}^6)}.
\end{aligned}$$
Hence, by Lemma \ref{lemma: basic properties of Wigner transform} $(ii)$, the two transforms are adjoint to each other;
\begin{equation}\label{eq: Weyl Wigner adjoint relation}
(\mathbf{Wig}^\hbar)^*=\mathbf{Weyl}^\hbar\quad\textup{and}\quad(\mathbf{Weyl}^\hbar)^*=\mathbf{Wig}^\hbar.
\end{equation}
\end{remark}

The two transformations are a well-known fundamental tool in semi-classical analysis. However, they do not preserve non-negativity, and it is not easy to prove boundedness from $L^p(\mathbb{R}^6)$ to $\mathfrak{L}_\hbar^p$ and vice versa except $p=2$. Alternatively, one can employ the \textit{Toeplitz quantization} (also called the Husimi quantization), defined by 
\begin{equation}\label{eq: Toeplitz quantization}
\mathbf{Toep}^\hbar[f]:=\frac{1}{(2\pi\hbar)^3}\iint_{\mathbb{R}^6}|\varphi_{(q,p)}^\hbar\rangle\langle\varphi_{(q,p)}^\hbar|f(q,p)dqdp
\end{equation}
for a distribution $f$, and the \textit{Husimi transform} defined by 
$$\mathbf{Hus}^\hbar[\gamma]:=\big\langle\varphi_{(q,p)}^\hbar\big|\gamma\big|\varphi_{(q,p)}^\hbar\big\rangle=\iint_{\mathbb{R}^6} \overline{\varphi_{(q,p)}^\hbar(x)}\gamma(x,x')\varphi_{(q,p)}^\hbar(x')dxdx'$$
for an operator $\gamma$, where $$\varphi_{(q,p)}^\hbar(x):=\frac{1}{(\pi\hbar)^{3/4}}e^{-\frac{|x-q|^2}{2\hbar}}e^{\frac{i(x-q)\cdot p}{\hbar}}$$ is a coherent state associated with $(q,p)\in\mathbb{R}^6$. Note that by definitions, these transforms preserve non-negativity, but they do not have inverse transforms. We also note that $\Tp$ and $\Hs$ are adjoint each other in the sense that 
\begin{equation}\label{eq:Husimi Toeplitz duality}
    (2\pi\hb)^3 \Tr(A^* \Tp[f])
    = \lg \Hs[A]|f \rg_{L^2_{q,p}(\R^6)}.
\end{equation}

The following lemma asserts that in the semi-classical limit $\hbar\to0$, the Weyl quantization (resp., the Wigner transform) gets closer to the Toeplitz quantization (resp., the Husimi transform).

\begin{lemma}[Regularization property of the Husimi transform and the Toeplitz quantization]\label{lemma: regularization of Weyl and Wigner transforms}
$\mathbf{Toep}^\hbar[f]=\mathbf{Weyl}^\hbar[G^\hbar*f]$ and  $\mathbf{Hus}^\hbar[\gamma]=G^\hbar*\mathbf{Wig}^\hbar[\gamma]$, where $G^\hbar(z)=\frac{e^{-|z|^2/\hbar}}{(\pi\hbar)^3}$ is a Gaussian function in $\mathbb{R}^6$.
\end{lemma}

\begin{remark}\label{remark: Wig Top}
$(i)$ By Lemma \ref{lemma: basic properties of Wigner transform} and Lemma \ref{lemma: regularization of Weyl and Wigner transforms}, the Toeplitz quantization and the Husimi transform are \textit{almost} inverses of each other in the sense that $\mathbf{Hus}^\hbar[\mathbf{Toep}^\hbar[f]]=G^\hbar*G^\hbar*f$. Similarly, for any $f\in L_{q,p}^2(\mathbb{R}^6)$, we have
$$\big\|\mathbf{Wig}^\hbar[\mathbf{Toep}^\hbar[f]]-f\big\|_{L_{q,p}^2(\mathbb{R}^6)}=\big\|G^\hbar*f-f\|_{L_{q,p}^2(\mathbb{R}^6)}\to 0.$$
$(ii)$ For a non-negative operator  $\gamma^\hbar\in\mathfrak{L}_\hbar^2$, if $\mathbf{Wig}^\hbar[\gamma^\hbar]\rightharpoonup f$ in $L_{q,p}^2(\mathbb{R}^6)$, then $f\geq 0$ (see \cite{Lions Paul 1993}). Indeed, $f$ is also the weak limit of non-negative $\mathbf{Hus}^\hbar[\gamma^\hbar]$, because by Lemma \ref{lemma: basic properties of Wigner transform} and Lemma \ref{lemma: regularization of Weyl and Wigner transforms}, 
$$\big\langle\mathbf{Wig}^\hbar[\gamma^\hbar]-\mathbf{Hus}^\hbar[\gamma^\hbar]\big|\varphi\big\rangle_{L_{q,p}^2(\mathbb{R}^6)}=(2\pi\hbar)^3\textup{Tr}\big(\gamma^\hbar \mathbf{Weyl}^\hbar[\varphi-G^\hbar*\varphi] \big)\to0.$$
Thus, $f$ is non-negative.
\end{remark}

\begin{proof}[Proof of Lemma \ref{lemma: regularization of Weyl and Wigner transforms}]
By the definition, the integral kernel of $\mathbf{Toep}^\hbar[f]$ can be written as 
$$\begin{aligned}
\mathbf{Toep}^\hbar[f](x,x')&=\frac{e^{-\frac{|x-x'|^2}{4\hbar}}}{(2\pi\hbar)^3}\iint_{\mathbb{R}^6} \frac{e^{-\frac{|q-\frac{x+x'}{2}|^2}{\hbar}}}{(\pi\hbar)^{3/2}}e^{\frac{i(x-x')\cdot p}{\hbar}}f(q,p)dqdp\\
&=\frac{e^{-\frac{|x-x'|^2}{4\hbar}}}{(2\pi\hbar)^3}\int_{\mathbb{R}^3} e^{\frac{i(x-x')\cdot p}{\hbar}}\bigg(\frac{e^{-\frac{|\cdot|^2}{\hbar}}}{(\pi\hbar)^{3/2}}*_qf\bigg)\bigg(\frac{x+x'}{2},p\bigg)dp\\
&=\frac{e^{-\frac{|x-x'|^2}{4\hbar}}}{\hbar^3}\bigg\{\mathcal{F}_p^{-1} \bigg(\frac{e^{-\frac{|\cdot|^2}{\hbar}}}{(\pi\hbar)^{3/2}}*_qf\bigg)(q,y)\bigg\}_{(q,y)=(\frac{x+x'}{2},\frac{x-x'}{\hbar})}.
\end{aligned}$$
Then, by the formula for the Fourier transform of Gaussian functions, it follows that
$$\begin{aligned}
\mathbf{Toep}^\hbar[f](x,x')&=\frac{1}{\hbar^3}\bigg\{e^{-\frac{\hbar|y|^2}{4}}\mathcal{F}_p^{-1} \bigg(\frac{e^{-\frac{|\cdot|^2}{\hbar}}}{(\pi\hbar)^{3/2}}*_qf\bigg)(q,y)\bigg\}_{(q,y)=(\frac{x+x'}{2},\frac{x-x'}{\hbar})}\\
&=\frac{1}{\hbar^3}\big(\mathcal{F}_p^{-1} (G^\hbar*f)(q,y)\big)\Big|_{(q,y)=(\frac{x+x'}{2},\frac{x-x'}{\hbar})}=\mathbf{Weyl}^\hbar[G^\hbar*f](x,x'),
\end{aligned}$$
since $\mathbf{Weyl}^\hbar[f](x,x')=\frac{1}{\hbar^3}(\mathcal{F}_p^{-1}f)(\frac{x+x'}{2},\frac{x-x'}{\hbar})$. Subsequently, we obtain $\mathbf{Hus}^\hbar[\gamma]=G^\hbar*\mathbf{Wig}^\hbar[\gamma]$ by duality \eqref{eq:Husimi Toeplitz duality}.
\end{proof}

Next, we prove the boundedness of the Toeplitz quantization and the Husimi quantization.
\begin{lemma}[Boundedness properties of the Toeplitz quantization/the Husimi quantization]\label{lem:Schatten Toeplitz}
For any $\al\in [1,\I)$ and $\si \in \mathbb{R}$, we have
\begin{align}
\big\|\sd^\si \mathbf{Toep}^\hbar[f] \sd^\si \big\|_{\FL^{\al}_\hb}
&\ls \|\lg p\rg^{2\si} f\|_{L^\al_{q,p}(\mathbb{R}^6)},\label{eq: weighted bound for Toeplitz (1)}\\
\big\|\lg \J_t^\hb\rg^{\si} \mathbf{Toep}^\hbar[f] \lg \J_t^\hb\rg^{\si} \big\|_{\FL^\al_\hb} 
&\ls \|\lg q-tp\rg^{2\si} f\|_{L^\al_{q,p}(\mathbb{R}^6)} \mbox{ for } |t|\le \frac{1}{\sqrt{\hb}}, \label{eq: weighted bound for Toeplitz (3)}\\
\big\|\langle p\rangle^{2\sigma}\mathbf{Hus}^\hbar[\gamma^\hbar]\big\|_{L_{q,p}^\alpha(\mathbb{R}^6)}
&\lesssim \|\langle \hbar\nabla\rangle^\sigma\gamma^\hbar\langle \hbar\nabla\rangle^\sigma\|_{\mathfrak{L}_\hbar^\alpha},\label{eq: weighted bound for Husimi (1)}\\
\big\|\langle q-tp\rangle^{2\sigma}\mathbf{Hus}^\hbar[\gamma^\hbar]\big\|_{L_{q,p}^\alpha(\mathbb{R}^6)}
&\lesssim \|\lg\J^\hb_t\rg^\sigma\gamma^\hbar\lg\J^\hb_t\rg^\sigma\|_{\mathfrak{L}_\hbar^\alpha}
\mbox{ for } |t|\le \frac{1}{\sqrt{\hb}}.\label{eq: weighted bound for Husimi (3)}
\end{align}
\end{lemma}

\begin{proof}
Note that by \eqref{eq:Husimi Toeplitz duality} and \eqref{eq: weighted bound for Toeplitz (1)}, for $\ps \in L^{\alpha'}_{q,p}$,
\begin{align*}
\big\lg \ps \big| \lg p \rg^{2\si} \Hs[\ga^\hb]\big\rg_{L^2_{q,p}}
&= (2\pi\hb)^3\Tr\Big(\sd^{-\si}\Tp[\lg p \rg^{2\si} \ps] \sd^{-\si} \sd^\si \ga^\hb \sd^\si\Big)\\
&\leq \big\|\sd^{-\si}\Tp[\lg p \rg^{2\si} \ps] \sd^{-\si}\big\|_{\mathfrak{L}_\hbar^{\alpha'}}\big\| \sd^\si \ga^\hb \sd^\si\big\|_{\mathfrak{L}_\hbar^\alpha}\\
&\leq \|\ps\|_{L_{q,p}^{\alpha'}}\big\| \sd^\si \ga^\hb \sd^\si\big\|_{\mathfrak{L}_\hbar^\alpha}.
\end{align*}
Hence, by duality, \eqref{eq: weighted bound for Husimi (1)} follows from \eqref{eq: weighted bound for Toeplitz (1)}. Similarly, \eqref{eq: weighted bound for Husimi (3)} is obtained from \eqref{eq: weighted bound for Toeplitz (3)}. Therefore, it is enough to show \eqref{eq: weighted bound for Toeplitz (1)} and \eqref{eq: weighted bound for Toeplitz (3)}.

For \eqref{eq: weighted bound for Toeplitz (1)}, we claim that 
\begin{equation}\label{eq: Toeplitz quantization bound 2}
\|\mathbf{Toep}^\hbar[f]\|_{\mathcal{B}(L^2)}\leq \|f\|_{L_{q,p}^\infty}.
\end{equation}
Indeed, for any $u\in L^2(\mathbb{R}^3)$, we have
$$\begin{aligned}
\big|\big\langle u\big| \mathbf{Toep}^\hbar[f]\big| u\big\rangle\big|&=\bigg|\frac{1}{(2\pi\hbar)^3}\iint_{\mathbb{R}^6}\iint_{\mathbb{R}^6}\overline{u(x)}u(x')\varphi_{(q,p)}^\hbar(x) \overline{\varphi_{(q,p)}^\hbar(x')}f(q,p)dxdx'dqdp\bigg|\\
&=\bigg|\frac{1}{(2\pi\hbar)^3}\iint_{\mathbb{R}^6} |(\varphi_{(0,p)}^\hbar*u)(q)|^2f(q,p)dqdp\bigg|\\
&\leq \|f\|_{L_{q,p}^\infty}\frac{1}{(2\pi\hbar)^3}\iint_{\mathbb{R}^6} |(\varphi_{(0,p)}^\hb*u)(q)|^2dqdp\\
&=\|f\|_{L_{q,p}^\infty}\iint_{\mathbb{R}^6} \bigg|\frac{1}{\hbar^{3/2}}\widehat{\varphi_{(0,p)}^\hbar}(\xi)\bigg|^2|\widehat{u}(\xi)|^2d\xi dp.
\end{aligned}$$
Note that 
\begin{equation}\label{eq: Fourier transform of coherent state}
\widehat{\varphi_{(0,p)}^\hbar}(\xi)=\int_{\mathbb{R}^3}\frac{1}{(\pi\hbar)^{3/4}}e^{-\frac{|x|^2}{2\hbar}}e^{\frac{ix\cdot p}{\hbar}}e^{-ix\cdot\xi}dx=(4\pi\hbar)^{\frac{3}{4}}e^{-\frac{\hbar}{2}|\xi-\frac{p}{\hbar}|^2}=(4\pi\hbar)^{\frac{3}{4}}e^{-\frac{|p-\hbar\xi|^2}{2\hbar}}.
\end{equation}
Hence, it follows that $|\langle u\big| \mathbf{Toep}^\hbar[f]\big| u\rangle|\leq\|f\|_{L_{q,p}^\infty}\|u\|_{L^2}^2$, which proves the claim.

On the other hand, for any $\sigma\in \R$, we have 
$$\begin{aligned}
\big\|\langle\hbar\nabla\rangle^\sigma\mathbf{Toep}^\hbar[f]\langle\hbar\nabla\rangle^\sigma\big\|_{\mathfrak{L}^1_\hb}&\leq\iint_{\mathbb{R}^6}\|\langle\hbar\nabla\rangle^\sigma\varphi_{(q,p)}^\hbar\|_{L_x^2}^2|f(q,p)|dqdp,\\
\big\|\langle \J_t^\hb \rangle^\sigma
\mathbf{Toep}^\hbar[f]
\langle \J^\hb_t\rangle^\sigma\big\|_{\mathfrak{L}^1_\hb}
&\leq\iint_{\mathbb{R}^6}\|\lg \J_t^\hb\rg^\sigma\varphi_{(q,p)}^\hbar\|_{L_x^2}^2|f(q,p)|dqdp.
\end{aligned}$$
By the Plancherel theorem with \eqref{eq: Fourier transform of coherent state}, we obtain that $$\|\langle\hbar\nabla\rangle^\sigma\varphi_{(q,p)}^\hbar\|_{L_x^2}^2=\frac{1}{(\pi\hbar)^{3/2}}\|\langle\xi\rangle^\sigma e^{-\frac{|p-\xi|^2}{2\hbar}}\|_{L_\xi^2}^2\sim \langle p\rangle^{2\sigma}.$$
Hence, we prove \eqref{eq: weighted bound for Toeplitz (1)} for $\alpha=1$.
Next, the simple calculation implies
\begin{equation*}
    \abs{\U^* \ph^\hb_{(q,p)} (x)}
    =\frac{1}{(\pi\hb)^{3/4}} \frac{1}{\bt^{3/2}} \exp\K{- \frac{|x-(q-tp)|^2}{2\hb(1+t^2)} }.
\end{equation*}
Hence, for $\si\ge 0$, we obtain
\begin{equation*}
\begin{aligned}
    &\|\lg\J^\hb_t\rg^{\si} \ph^\hb_{(q,p)}\|_{L^2_x}
    = \frac{1}{(\pi\hb)^{3/4}} \frac{1}{\bt^{3/2}} 
    \No{\bx^{\si} \exp\K{-\frac{|x-(q-tp)|^2}{2\hb(1+t^2)} } }_{L^2_x} \\
    &\quad =\frac{1}{(\pi\hb)^{3/4}} \frac{1}{\bt^{3/2}} 
    \No{\lg x+(q-tp) \rg^{\si} \exp\K{ -\frac{|x|^2}{2\hb(1+t^2)} } }_{L^2_x} \\
    &\quad \ls \lg q - tp \rg^{\si} \No{\lg \sqrt{2\hb}\bt x \rg^{\si} e^{-|x|^2}}_{L^2_x}
    \ls\lg q - tp \rg^{\si}.
\end{aligned}
\end{equation*}
if $|t| \le \frac{1}{\sqrt{\hb}}$, where we used $\lg a + b \rg \ls \lg a \rg \lg b \rg$.
When $\si \le 0$, we can do the same calculation as above by using $\lg a + b \rg^{-1} \ls \lg a \rg \lg b \rg^{-1}$.
Hence, we have \eqref{eq: weighted bound for Toeplitz (1)} and \eqref{eq: weighted bound for Toeplitz (3)} for $\al=1$.
Therefore, by the claim and the complex interpolation with \eqref{eq: Toeplitz quantization bound 2}, we complete the proof of the lemma.
\end{proof}

\subsection{Wigner transform of the linear Schr\"odinger equation}
The Wigner transform associates the quantum linear flow with the transport equation. 
\begin{lemma}[Wigner transform of the linear Schr\"odinger flow]\label{lem:equation for f hbar}
For $\gamma_0^\hbar\in\mathfrak{L}_\hbar^2$, let $\gamma^\hbar(t)=\mathcal{U}_{V^\hbar}^\hbar(t)\gamma_0^\hbar \mathcal{U}_{V^\hbar}^\hbar(t)^*$. Then, $f^\hbar(t)=\mathbf{Wig}^\hbar[\gamma^\hbar(t)] \in C(\R;L^2_{q,p}(\mathbb{R}^6))$ satisfies
	\begin{equation}\label{eq: reformulated LS}
		\left\{\begin{aligned}
			\partial_t f^\hbar+p\cdot\nabla_q f^\hbar-\nabla V^\hbar\cdot\nabla_p f^\hbar&=-i\mathcal{R}^\hbar[\gamma^\hbar; V^\hbar],\\
			f^\hbar(0)&=\mathbf{Wig}^\hbar[\gamma_0^\hbar]
		\end{aligned}\right.
	\end{equation}
	in the sense of the mild solution in Remark \ref{remark: mild solution definition of the Vlasov equation}, where the remainder term is given by 
	\begin{equation}\label{eq: reformulated LS remainder}
		\begin{aligned}
		\mathcal{R}^\hbar[\gamma; V](q,p)&:=\mathcal{F}_{y\to p}\Bigg[\hbar^2\mathcal{Q}^\hbar[V](q,y)\gamma\bigg(q+\frac{\hbar y}{2},q-\frac{\hbar y}{2}\bigg)\Bigg]
		\end{aligned}
	\end{equation}
    and
    $$\mathcal{Q}^\hbar[V](q,y):=V\bigg(q+\frac{\hbar y}{2}\bigg)-V\bigg(q-\frac{\hbar y}{2}\bigg)-\nabla V(q)\cdot \hbar y.$$
\end{lemma}

\begin{remark}\label{remark: mild solution definition of the Vlasov equation}

Similarly to the definition of a global solution to the Vlasov equation (see Remark \ref{definition: scattering solutions to Vlasov}), we say that $f^\hbar(t)$ is an $H^{-1}$-global solution to  \eqref{eq: reformulated LS} with initial data $f^\hbar(0)$ if $g^\hbar(t)=\mathcal{U}(-t)f^\hbar(t)$ satisfies 
$$\begin{aligned}
g^\hbar(t)&=f^\hbar(0) + \int_0^t  \CU(-t_1) \Big\{(\na_q V^\hb) \cdot \nabla_pf^\hbar-i\mathcal{R}^\hbar[\gamma^\hbar; V^\hbar]\Big\}(t_1) dt_1\\
&\in C(\mathbb{R}_t;H_{q,p}^{-1}(\mathbb{R}^6)),
\end{aligned}$$
where $\mathcal{U}(t)$ is the free transport flow (see \eqref{eq: free transport flow, definition}). 
\end{remark}

In \eqref{eq: reformulated LS},  the remainder is small due to the following bounds.

\begin{lemma}[Remainder estimates]\label{lem: remainder estimate} 
For $0\leq s\leq 1$ and $2\leq r<\infty$, we have
\begin{align}
\big\|\langle \nabla_p\rangle^{-1}\mathcal{R}^\hbar[\gamma; V]\big\|_{L_{q,p}^2(\mathbb{R}^6)}&\lesssim\|\nabla V\|_{L^\infty}\|\gamma\|_{\mathfrak{L}_\hbar^2},\label{eq: remainder estimate 1}\\
\big\|\langle \nabla_q\rangle^{-2}\langle \nabla_p\rangle^{-2}\mathcal{R}^\hbar[\gamma; V]\big\|_{L_{q,p}^2(\mathbb{R}^6)}&\lesssim \hbar^s\|\nabla V\|_{\dot{W}_q^{s,r}}\|\gamma\|_{\mathfrak{L}_\hbar^2}.\label{eq: remainder estimate 2}
\end{align}
\end{lemma}

\begin{proof}
By the fundamental theorem of calculus, $|\mathcal{Q}^\hbar[V](q,y)|\leq 2\hbar\|\nabla V\|_{L^\infty}|y|$ holds. Thus, by the Plancherel theorem, we prove that 
$$\begin{aligned}
\big\|\langle \nabla_p\rangle^{-1}\mathcal{R}^\hbar[\gamma; V]\big\|_{L_{q,p}^2}&\sim\bigg\|\frac{\hbar^2}{\langle y\rangle}\mathcal{Q}^\hbar[V](q,y)\gamma\bigg(q+\frac{\hbar y}{2},q-\frac{\hbar y}{2}\bigg)\bigg\|_{L_{q,y}^2}\\
&\lesssim\hbar^3 \|\nabla V\|_{L^\infty}\bigg\|\gamma\bigg(q+\frac{\hbar y}{2},q-\frac{\hbar y}{2}\bigg)\bigg\|_{L_{q,y}^2}\\
&\sim\|\nabla V\|_{L^\infty}(2\pi\hbar)^\frac{3}{2}\|\gamma\|_{\mathfrak{S}^2}=\|\nabla V\|_{L^\infty}\|\gamma\|_{\mathfrak{L}_\hbar^2}.
\end{aligned}$$
For \eqref{eq: remainder estimate 2}, similarly but using the Sobolev inequality, we obtain
$$\begin{aligned}
\big\|\langle \nabla_q\rangle^{-2}\langle \nabla_p\rangle^{-2}\mathcal{R}^\hbar[\gamma; V]\big\|_{L_{q,p}^2}&\sim\bigg\|\langle \nabla_q\rangle^{-2}\bigg(\frac{\hbar^2}{\langle y\rangle^2}\mathcal{Q}^\hbar[V](q,y)\gamma\bigg(q+\frac{\hbar y}{2},q-\frac{\hbar y}{2}\bigg)\bigg)\bigg\|_{L_{q,y}^2}\\
&\lesssim\bigg\|\frac{\hbar^2}{\langle y\rangle^2}\mathcal{Q}^\hbar[V](q,y)\gamma\bigg(q+\frac{\hbar y}{2},q-\frac{\hbar y}{2}\bigg)\bigg\|_{L_y^2L_q^{\frac{2r}{r+2}}}\\
&\leq\hbar^2\bigg\|\frac{1}{\langle y\rangle^2}\mathcal{Q}^\hbar[V](q,y)\bigg\|_{L_y^\infty L_q^r}\bigg\|\gamma\bigg(q+\frac{\hbar y}{2},q-\frac{\hbar y}{2}\bigg)\bigg\|_{L_y^2L_q^2}.
\end{aligned}$$
Then, we observe that by the fundamental theorem of calculus, $\mathcal{Q}^\hbar[V]$ can be written more precisely as 
$$\begin{aligned}
\mathcal{Q}^\hbar[V](q,y)&=\int_{-\frac{1}{2}}^{\frac{1}{2}} \frac{d}{d\th_1}\big(V(q+\hbar y\th_1)\big)d\th_1-\nabla V(q)\cdot \hbar y\\
&=\bigg\{\int_{-\frac{1}{2}}^{\frac{1}{2}} \K{(\nabla V)(q+\hbar y\th_1)-\nabla V(q)} d\th_1\bigg\}\cdot \hbar y\\
&=\hbar y\bigg\{\int_{-\frac{1}{2}}^{\frac{1}{2}}\int_0^{\theta_1}(\nabla^2 V)(q+\hbar y\th_2)d\theta_2 d\th_1\bigg\} \hbar y.
\end{aligned}$$
Hence, it follows that $\|\mathcal{Q}^\hbar[V](q,y)\|_{L_q^r}\lesssim \hbar|y|\|V\|_{\dot{W}_q^{1,r}}$ and $\|\mathcal{Q}^\hbar[V](q,y)\|_{L_q^r}\lesssim (\hbar|y|)^2\|V\|_{\dot{W}_q^{2,r}}$, and consequently, interpolating them, we obtain 
$$\|\mathcal{Q}^\hbar[V](q,y)\|_{L_q^r}\lesssim(\hbar|y|)^{1+s}\|V\|_{\dot{W}_q^{1+s,r}}\lesssim(\hbar|y|)^{1+s}\|\nabla V\|_{\dot{W}_q^{s,r}}$$
for $0\leq s\leq 1$. Therefore, we conclude that
$$\big\|\langle \nabla_q\rangle^{-2}\langle \nabla_p\rangle^{-2}\mathcal{R}^\hbar[\gamma; V]\big\|_{L_{q,p}^2}\lesssim\hbar^s\|\nabla V\|_{\dot{W}_q^{s,r}}\|\gamma\|_{\mathfrak{L}_\hbar^2},$$
where we used that $\hbar^3\|\gamma(q+\frac{\hbar y}{2},q-\frac{\hbar y}{2})\|_{L_{q,y}^2}\sim\|\gamma\|_{\mathfrak{L}_\hbar^2}$.
\end{proof}

In the free case $V^\hbar\equiv0$, \eqref{eq: reformulated LS} is the free transport equation, and
$\mathbf{Wig}^\hbar[\CU^\hb(t)\gamma_0^\hbar\CU^\hb(t)^*]=\E \mathbf{Wig}^\hbar[\gamma_0^\hbar]$. The following lemma asserts that a similar identity holds for time-dependent operators.

\begin{lemma}[Wigner transform of the backward free flow]\label{lemma: Wigner transform of free evolution}
	For any $\ga(t)\in C(\R;\FS^2)$, we have
$$\mathbf{Wig}^\hbar[\Uinv \gamma(t) \U ]
=\CU(t)^* \mathbf{Wig}^\hbar[\gamma(t)].$$
\end{lemma}

\begin{proof}
By the integral representation of the free linear Schr\"odinger flow and collecting the terms with the $y$-variable, we write 
$$\begin{aligned}
&\mathbf{Wig}^\hbar[\Uinv \gamma(t) \U ](q,p)\\
&=\iiint_{\mathbb{R}^{9}}\frac{1}{(2\pi\hbar|t|)^3}e^{\frac{|q+\frac{y}{2}-z|^2-|q-\frac{y}{2}-z'|^2}{2i\hbar t}}\gamma(t,z,z') e^{-\frac{iy\cdot p}{\hbar}}dydzdz'\\
&=\iiint_{\mathbb{R}^{9}}\frac{1}{(2\pi\hbar|t|)^3}e^{\frac{|q-z|^2-|q-z'|^2}{2i\hbar t}}\gamma(t,z,z') e^{-\frac{iy}{2t\hbar}\cdot (2q+2tp-z-z')}dydzdz'.
\end{aligned}$$
Integrating out the $y$- and then the $z'$-variables, we obtain 
$$\begin{aligned}
\mathbf{Wig}^\hbar[\Uinv \gamma(t) \U ](q,p)
&=8 \int_{\mathbb{R}^3}e^{\frac{|q-z|^2-|q+2tp-z|^2}{2i\hbar t}}\gamma(t,z,2q+2tp-z) dz.
\end{aligned}$$
Then, changing the variable $z=q+tp+\frac{y}{2}$, we obtain that 
$$\begin{aligned}
&\mathbf{Wig}^\hbar[\Uinv \gamma(t) \U ](q,p)\\
&=\int_{\mathbb{R}^3}e^{\frac{|tp+\frac{y}{2}|^2-|tp-\frac{y}{2}|^2}{2i\hbar t}}\gamma\bigg(t,q+tp+\frac{y}{2},q+tp-\frac{y}{2}\bigg) dy\\
&=\int_{\mathbb{R}^3}\gamma\bigg(t,(q+tp)+\frac{y}{2},q+tp-\frac{y}{2}\bigg) e^{-\frac{iy\cdot p}{\hbar}}dy
=\CU(t)^* \mathbf{Wig}^\hbar[\gamma(t) ](q,p).
\end{aligned}$$
\end{proof}

\section{Semi-classical limit of scattering states: Proof of Theorem \ref{thm: semiclassical limit}}\label{sec: semi-classical limit of scattering states}

This section is devoted to our second main result (Theorem \ref{thm: semiclassical limit}) whose proof is broken down into three steps. First, having small-data global quantum states in Theorem \ref{thm: uniform bounds for nonlinear solutions}, in Section \ref{sec: compactness}, we prove the compactness of their Wigner transforms and a priori bounds for their limits (Proposition \ref{prop: compactness of fh}). Next, in Section \ref{sec:construction of f}, we show that the limits are classical scattering states for the Vlasov equation (Proposition \ref{prop:Vlasov sol}). Finally, in Section \ref{sec: correspondence of the scattering states}, we establish the weak convergence from quantum to classical scattering states via the Wigner transform (Proposition \ref{prop:correspondence of the scattering state}). In Section \ref{sec:scattering Vlasov}, as a byproduct, we present a new scattering result for the Vlasov equation  (Theorem \ref{thm:scattering Vlasov}).

\subsection{Compactness for quantum states and a priori bounds for their limits: Proof of Theorem \ref{thm: semiclassical limit} \texorpdfstring{$(i)$}{(i)}}\label{sec: compactness}

To begin with, we show the compactness of Wigner transforms of small scattering quantum states, and deduce several useful bounds of their limits by the Banach-Alaoglu theorem.

\begin{proposition}[Compactness of the Wigner transforms of small-data solutions to NLH]\label{prop: compactness of fh}
For $\hbar\in(0,1]$, let $$g^\hbar(t):=\mathcal{U}(-t)\mathbf{Wig}^\hbar[\gamma^\hbar(t)],$$ where $\gamma^\hbar(t)$ is a small-data global solution to NLH \eqref{eq: NLH} constructed in Theorem \ref{thm: uniform bounds for nonlinear solutions}. Then, the following hold.
\begin{enumerate}[$(i)$]
\item On the whole time interval $\mathbb{R}$, $\mathcal{G}:=\{g^\hbar(t)\}_{\hbar\in(0,1]}$ is uniformly bounded in $L^2_{q,p}(\mathbb{R}^6)$ and equicontinuous with respect to the metric $\mathbf{d}^w$ (see \eqref{eq: weak-topology metric}). As a consequence, for any $\{\hbar_j\}_{j=1}^\infty$ such that $\hbar_j\to0$, there exists $g(t)\in C_t(\mathbb{R};  w-B_{\eta_0})$ such that passing to a subsequence, $g^{\hbar_j}(t)$ converges to $g(t)$ in $ C_t(I; w-B_{\eta_0})$ for any finite interval $I$.
\item The limit $g(t)$, obtained in $(i)$, is non-negative. Moreover, $f(t):=\mathcal{U}(t)g(t)$ obeys
\begin{equation}\label{eq: weighted norm bound for the weak limit f}
\begin{aligned}
&\|f(t)\|_{ L^\I_t(\mathbb{R};L_{q,p}^1(\mathbb{R}^6))}+\|\langle p\rangle^{1+a+2\epsilon}
f(t)\|_{ L^\I_t(\mathbb{R};L_{q,p}^{\frac{3}{2-a-\epsilon}}(\mathbb{R}^6))}\\
&\qquad \qquad +\|\langle q-tp\rangle^{1+a+2\epsilon} f(t)\|_{ L^\I_t(\mathbb{R};L_{q,p}^{\frac{3}{2-a-\epsilon}}(\mathbb{R}^6))}\lesssim \eta_0
\end{aligned}
\end{equation}
and
\begin{equation}\label{eq: decay of limit of densities}
\sup_{t\in\mathbb{R}}\langle t\rangle^{1+a}\|\nabla w*\rho_{f(t)}\|_{L^\infty_x(\mathbb{R}^3)}+ \sup_{t\in\mathbb{R}}\langle t\rangle^{1+a+\epsilon}\|\rho_{f(t)}\|_{L^{\frac{3}{2-a-\epsilon}}_x}\lesssim \eta_0.
\end{equation}
\end{enumerate}
\end{proposition}

\begin{proof}
For $(i)$, it is easy to see that the set $\mathcal{G}$ is uniformly bounded on $\mathbb{R}$ in $L_{q,p}^2$, because by the isometry of the Wigner transform (Lemma \ref{lemma: basic properties of Wigner transform}) and the conservation of the $\mathfrak{L}_h^2$-norm for $\gamma^\hbar(t)$,  
\begin{equation}\label{fh L^2 norm}
\|g^\hbar(t)\|_{L_{q,p}^2}=\|f^\hbar(t)\|_{L_{q,p}^2}=\|\gamma^\hbar(t)\|_{\mathfrak{L}_h^2}=\|\gamma_0^\hbar\|_{\mathfrak{L}_h^2}\lesssim \eta_0.
\end{equation}
For equicontinuity, assuming $t_2\geq t_1$ without loss of generality, we consider the difference 
$$\mathbf{d}^w(g^\hbar(t_2),g^\hbar(t_1))=\sum_{n=1}^\infty 2^{-n}\abs{\big\langle \psi_n\big|g^\hbar(t_2)-g^\hbar(t_1)\big\rangle_{L_{q,p}^2(\mathbb{R}^6)}}.$$
Given any small number $\delta>0$, for each $n\geq1$, we choose $\psi_n^\delta\in C_c^\infty(\mathbb{R}^6)$ such that $\|\psi_n^\delta-\psi_n\|_{L_{q,p}^2}\leq\delta$. Then, by the uniform bound \eqref{fh L^2 norm}, one can find $N=N_\delta\gg1$ such that 
$$\mathbf{d}^w(g^\hbar(t_2),g^\hbar(t_1))
=\sum_{n=1}^N 2^{-n}\abs{\big\langle \psi_n^\delta\big|g^\hbar(t_2)-g^\hbar(t_1)\big\rangle_{L_{q,p}^2(\mathbb{R}^6)}}+O(\delta).$$
Hence, it suffices to show that as a function of $t$, $\langle \psi|g^\hbar(t)\rangle_{L_{q,p}^2}$ is uniformly equicontinuous on $\mathbb{R}$, assuming that $\psi\in C_c^\infty(\mathbb{R}^6)$. For the proof, by the Duhamel representation for $f^\hb(t)$ (see Lemma \ref{lem:equation for f hbar}), we write 
$$\begin{aligned}
\big\langle \psi\big|g^\hbar(t_2)-g^\hbar(t_1)\big\rangle_{L_{q,p}^2}&=-\int_{t_1}^{t_2}\big\langle \nabla_p \CU(t')\psi\big|(\nabla \Phi^\hbar)(t')f^\hbar(t')\big\rangle_{L_{q,p}^2} dt'\\
&\quad-i\int_{t_1}^{t_2} \big\langle \CU(t')\psi\big| \CR^\hbar[\ga^\hb;\Phi^\hbar](t')\big\rangle_{L_{q,p}^2} dt'\\
&=:\textup{(I)}+\textup{(II)}.
\end{aligned}$$
For $\textup{(I)}$, we observe that $(\nabla_p\CU(t')\psi)(q,p)=-t'(\nabla_q\psi)(q-t'p,p)+(\nabla_p\psi)(q-t'p,p)$. Hence, it follows from the uniform decay estimate \eqref{eq: dispersion estimate'} for the potential $\Phi^\hbar$, the isometry of the Wigner transform (Lemma \ref{lemma: basic properties of Wigner transform}) for $f^\hbar(t)=\mathbf{Wig}^\hbar[\gamma^\hbar(t)]$ and the conservation law bound \eqref{fh L^2 norm}  that 
$$\begin{aligned}
|\textup{(I)}|&\leq\int_{t_1}^{t_2}\|\nabla_p\CU(t')\psi\|_{L_{q,p}^2}\|\nabla \Phi^\hbar(t')\|_{L_q^\infty}\|f^\hbar(t')\|_{L_{q,p}^2} dt'\\
&\lesssim\int_{t_1}^{t_2}\langle t'\rangle\|\psi\|_{H_{q,p}^1}\frac{\eta_0}{\langle t'\rangle^{1+a}}\eta_0dt'\lesssim \eta_0^2 (t_2-t_1)\|\psi\|_{H_{q,p}^1}.
\end{aligned}$$
On the other hand, for $\textup{(II)}$, applying \eqref{eq: remainder estimate 1} to $\CR^\hbar[\ga^\hb;\Phi^\hbar](t')$ and estimating as before, we obtain
$$\begin{aligned}
|\textup{(II)}|&\leq\int_{t_1}^{t_2} \|\langle\nabla_p\rangle\CU(t')\psi\|_{L_{q,p}^2} \|\langle\nabla_p\rangle^{-1}\CR^\hbar[\ga^\hb;\Phi^\hbar](t')\|_{L_{q,p}^2} dt'\\
&\lesssim\int_{t_1}^{t_2}\langle t'\rangle\|\psi\|_{H_{q,p}^1}\frac{\eta_0}{\langle t'\rangle^{1+a}}\eta_0dt'\lesssim \eta_0^2 (t_2-t_1)\|\psi\|_{H_{q,p}^1}.
\end{aligned}$$
Therefore, we prove the desired uniform equicontinuity for $\langle \psi|g^\hbar(t)\rangle_{L_{q,p}^2}$ as well as that for $\mathcal{G}=\{g^\hbar(t)\}_{\hbar\in(0,1]}$.

Now, by the Arzelà–Ascoli theorem, we note that any sequence $\{\hbar_j\}_{j=1}^\infty$ with $\hbar_j\to0$ has a subsequence $\{\hbar_{1;j}\}_{j=1}^\infty$ such that $g^{\hbar_{1;j}}(t)\to g(t)$ in $C_t([-1,1]; w-B_{\eta_0})$. Then, we also note that for an integer $k\geq 1$, if $g^{\hbar_{k;j}}(t)\to g(t)$ in $C_t([-k,k]; w-B_{\et})$ for some sequence $\{\hbar_{k;j}\}_{j=1}^\infty$, then $\{\hbar_{k;j}\}_{j=1}^\infty$ has a subsequence $\{\hbar_{k+1;j}\}_{j=1}^\infty$ such that $g^{\hbar_{k+1;j}}(t)\to g(t)$ in $C_t([-(k+1),k+1]; w-B_{\et})$. Therefore, by the standard diagonal argument, we conclude that $g^{\hbar_{j;j}}(t)\to g(t)$ in $C_t([-T,T]; w-B_{\et})$ for any $T>0$ for some $g\in C_t(\mathbb{R}; w-B_{\et})$. For convenience, we denote $\hbar_j=\hbar_{j;j}$ with abuse of notation.

For $(ii)$, we observe that for each $t$, $g(t)$ coincides with the weak limit of $$\mathcal{U}(-t)\mathbf{Hus}^{\hbar_j}[\gamma^{\hbar_j}(t)]=\mathcal{U}(-t)\big(G^{\hbar_j}*\mathbf{Wig}^{\hbar_j}[\gamma^{\hbar_j}(t)]\big)$$
in $L_{q,p}^2$ (see Lemma \ref{lemma: regularization of Weyl and Wigner transforms}, Lemma \ref{lemma: Wigner transform of free evolution} and Remark \ref{remark: Wig Top} $(ii)$), and thus $g(t,q,p)\geq0$.
We also note that for each $t$, by the second inequality in \eqref{eq: weighted norm uniform bound} and Lemma \ref{lem:Schatten Toeplitz},
\begin{equation}\label{eq: Husimi of gamma weighted norm bound}
\begin{aligned}
\big\|\langle p\rangle^{1+a+2\epsilon}\mathcal{U}(-t)\mathbf{Hus}^{\hbar_j}[\gamma^{\hbar_j}(t)]\big\|_{L_{q,p}^{\frac{3}{2-a-\epsilon}}}
&=\big\|\langle p\rangle^{1+a+2\epsilon}\mathbf{Hus}^{\hbar_j}[\gamma^{\hbar_j}(t)]\big\|_{L_{q,p}^{\frac{3}{2-a-\epsilon}}},\\
\big\|\langle q\rangle^{1+a+2\epsilon}\mathcal{U}(-t)\mathbf{Hus}^{\hbar_j}[\gamma^{\hbar_j}(t)]\big\|_{L_{q,p}^{\frac{3}{2-a-\epsilon}}}
&=\big\|\langle q-tp\rangle^{1+a+2\epsilon}\mathbf{Hus}^{\hbar_j}[\gamma^{\hbar_j}(t)]\big\|_{L_{q,p}^{\frac{3}{2-a-\epsilon}}},
\end{aligned}
\end{equation}
are bounded uniformly for all $|t|\leq\frac{1}{\sqrt{\hbar}}$. Thus, it follows from the Banach-Alaoglu theorem that 
$$\|g(t)\|_{L^\I_t(\mathbb{R};L_{q,p}^1)}
+\|\langle p\rangle^{1+a+2\epsilon} g(t)\|_{L^\I_t(\mathbb{R};L_{q,p}^{\frac{3}{2-a-\epsilon}})}+\|\langle q\rangle^{1+a+2\epsilon} g(t)\|_{L^\I_t(\mathbb{R};L_{q,p}^{\frac{3}{2-a-\epsilon}})}\lesssim \eta_0,$$
which deduces \eqref{eq: weighted norm bound for the weak limit f}.

For \eqref{eq: decay of limit of densities}, we define $\rho_j:=\rho_{f^{\hbar_j}}=\rh^{\hb_j}_{\ga^{\hb_j}}$. Then, interpolating the inequalities in \eqref{eq: dispersion estimate}, we have $\langle t\rangle^{3/2}\|\rho_j(t)\|_{L^2_x}\lesssim \eta_0$ and thus, $\|\rho_j(t)\|_{L_t^2(\mathbb{R};L^2_x)}\lesssim \eta_0$. Hence, extracting a subsequence, we obtain $\rho_j\rightharpoonup\rho_\infty$ in  $L_t^2(\mathbb{R};L^2_x)$ for some $\rho_\infty$. Moreover, by the Banach-Alaoglu theorem, the limit satisfies $$\sup_{t\in\mathbb{R}}\langle t\rangle^{1+a+\epsilon}\|\rho_\infty(t)\|_{L^{\frac{3}{2-a-\epsilon}}_x}\lesssim \eta_0.$$
Therefore, it suffices to show that $\rho_\infty(t)=\rho_{f}(t)$ in $L_x^2$ for each $t\in [-T,T]$ with any finite $T>0$. Indeed, by the uniqueness of distributions, it is enough to prove that for any $\psi\in C_c^\infty(\mathbb{R}^3)$,
$$\langle \psi|\rho_j(t)-\rho_{f(t)}\rangle_{L_x^2}\to 0.$$
Note that by interpolating the bounds in \eqref{eq: weighted norm bound for the weak limit f}, we have
$$\|\langle p\rangle^{\frac{3(1+a+2\epsilon)}{2(1+a+\epsilon)}} f(t)\|_{C_t([-T,T]; L_{q,p}^2)}\lesssim \eta_0$$
so that by H\"older's inequality, 
$$\sup_{t\in[-T,T]}\|\rho_{f(t)}-\rho_{\chi_R(p)f(t)}\|_{L_q^2}=o_R(1),$$
where $\chi_R(p)$ is a smooth high momentum cut-off such that $\chi_R\equiv1$ in $|p|\leq R$ but $\chi_R\equiv0$ in $|p|\geq 2R$, and $o_R(1)$ is a small number converging to zero as $R\to\infty$. On the other hand, by Lemma \ref{lemma: regularization of Weyl and Wigner transforms}, we have 
$$\begin{aligned}
\big|\big\langle \psi\big|\rho_j(t)-\rho_{\chi_R(p)f^{\hbar_j}(t)}\big\rangle_{L_x^2}\big|&\leq \|\psi\|_{L_x^2}\|\rho_{\mathbf{Hus}^{\hbar_j}[\gamma^{\hbar_j}(t)]}-\rho_{\chi_R(p)\mathbf{Hus}^{\hbar_j}[\gamma^{\hbar_j}(t)]}\|_{L_x^2}\\
&\quad+\big|\big\langle \psi\big|\rho_{f^{\hbar_j}(t)}-\rho_{G^{\hbar_j}*f^{\hbar_j}(t)}\big\rangle_{L_x^2}\big|\\
&\quad+\big|\big\langle \psi\big|\rho_{\chi_R(p)f^{\hbar_j}(t)}-\rho_{\chi_R(p)G^{\hbar_j}*f^{\hbar_j}(t)}\big\rangle_{L_x^2}\big|.
\end{aligned}$$
Since $\rho_{G^{\hbar}*f^{\hbar}(t)}=\frac{e^{-|\cdot|^2/\hbar}}{(\pi\hbar)^{3/2}}*\rho_{f^{\hbar}}(t)$ and the convolution with $\frac{e^{-|\cdot|^2/\hbar}}{(\pi\hbar)^{3/2}}$ is an approximation to identity, we obtain
$$\begin{aligned}
\big\langle \psi\big|\rho_{f^{\hbar_j}(t)}-\rho_{G^{\hbar_j}*f^{\hbar_j}(t)}\big\rangle_{L_x^2}&=\big\langle \psi-\tfrac{e^{-|\cdot|^2/\hbar}}{(\pi\hbar)^{3/2}}*\psi\big|\rho_{f^{\hbar_j}(t)}\big\rangle_{L_x^2}=o_j(1),
\end{aligned}$$
where $o_j(1)$ is a small number that converges to zero as $j\to\infty$. Similarly, one can show that
$$\begin{aligned}
&\big\langle \psi\big|\rho_{\chi_R(p)f^{\hbar_j}(t)}-\rho_{\chi_R(p)G^{\hbar_j}*f^{\hbar_j}(t)}\big\rangle_{L_x^2}\\
&=\Big\langle\psi(q)\chi_R(p)-G^{\hbar_j}*\big(\psi(q)\chi_R(p)\big)\big|f^{\hbar_j}(t)\Big\rangle_{L_{q,p}^2}=o_j(1).
\end{aligned}$$ Subsequently, using the weighted norm bound \eqref{eq: Husimi of gamma weighted norm bound}, we prove that 
$$\begin{aligned}
\big|\langle \psi|\rho_j(t)-\rho_{\chi_R(p)f^{\hbar_j}(t)}\rangle_{L_x^2} \big|
&\le\|\psi\|_{L^2_x}\|\rho_{\mathbf{Hus}^{\hbar_j}[\gamma^{\hbar_j}(t)]}-\rho_{\chi_R(p)\mathbf{Hus}^{\hbar_j}[\gamma^{\hbar_j}(t)]}\|_{L_q^2}+o_j(1)\\
&=o_j(1)+o_R(1).
\end{aligned}$$
Therefore, by the weak convergence of $g^{\hbar_j}(t)$ in $(i)$, we obtain
$$\begin{aligned}
\langle \psi|\rho_j(t)-\rho_{f(t)}\rangle_{L_x^2}&=\langle \psi|\rho_{\chi_R(p)(f^{\hbar_j}-f)(t)}\rangle_{L_x^2}+o_j(1)+o_R(1)\\
&=\big\langle \psi(q+tp)\chi_R(p)\big|(g^{\hbar_j}-g)(t)\big\rangle_{L_{q,p}^2}+o_j(1)+o_R(1)\\
&\to o_R(1)
\end{aligned}$$
as $j\to\infty$. Since $R\gg 1$ can be arbitrary, by the uniqueness of the weak limit, we conclude that $\rho_\infty(t)=\rho_f(t)$ in $L_x^2$.
\end{proof}

\subsection{Construction of a global scattering solution to the Vlasov equation: Proof of Theorem \ref{thm: semiclassical limit} \texorpdfstring{$(ii)$}{(ii)}}\label{sec:construction of f}

In the previous subsection, it is shown that $f=\mathcal{U}(t)g(t)$ obeys the bounds in Proposition \ref{prop: compactness of fh}. For Theorem \ref{thm: semiclassical limit} $(ii)$, we show the following proposition. 

\begin{proposition}[Derivation of classical scattering states]\label{prop:Vlasov sol}
Let $g(t)\in C_t(\mathbb{R}; w-B_{\et})$ be the limit function constructed in Proposition \ref{prop: compactness of fh}, and $f=\mathcal{U}(t)g(t)$. Then, $g(t)$ is a solution to the Vlasov equation in a moving reference frame 
\begin{equation}\label{eq: Vlasov-moving frame'}
\begin{aligned}
g(t)&=f_0+\int_0^t \CU(-t_1) \Big\{\nabla_q (w*\rho_f)\cdot\nabla_pf\Big\}(t_1)dt_1 \\
&=f_0+\int_0^t (-t_1\na_q+\na_p)\CU(-t_1) \Big\{\nabla_q (w*\rho_f)\cdot f\Big\}(t_1)dt_1
\end{aligned}
\end{equation}
in $C_t(\mathbb{R}; H_{q,p}^{-1}(\mathbb{R}^6))$. Moreover, there exists $f_\pm\in H_{q,p}^{-1}(\mathbb{R}^6)$ such that
$$\lim_{t\to\pm\infty}\|\mathcal{U}(-t)f(t)-f_\pm\|_{H_{q,p}^{-1}(\mathbb{R}^6)}=0.$$
\end{proposition}

For the proof, the following elementary lemmas will be useful. 

\begin{lemma}\label{lemma: integral estimates for derivation of VP}
If $1<a<2$ and $\epsilon>0$ is sufficiently small, then 
$$\begin{aligned}
\big|\big\langle \psi\big|(\nabla_q w \ast \zeta) f\big\rangle_{L_{q,p}^2(\mathbb{R}^6)}\big|&\leq \|\psi\|_{L_{q,p}^2(\mathbb{R}^6)}\|\nabla_q w\|_{L_q^\infty + L_q^{\frac{3}{1+a+\epsilon}}(\mathbb{R}^3)} \|\zeta\|_{L_q^1\cap L_q^{\frac{3}{2-a-\epsilon}}(\mathbb{R}^3)}\|f\|_{L_{q,p}^2(\mathbb{R}^6)},\\
\big|\big\langle \psi\big|\CR^\hb[\ga;w \ast \zeta]\big\rangle_{L_{q,p}^2}\big|&\lesssim\|\langle\nabla_p\rangle \psi\|_{L_{q,p}^2(\mathbb{R}^6)}\|\nabla_q w\|_{L_q^\infty + L_q^{\frac{3}{1+a+\epsilon}}(\mathbb{R}^3)} \|\zeta\|_{L_q^1\cap L_q^{\frac{3}{2-a-\epsilon}}(\mathbb{R}^3)}\|\gamma\|_{\mathfrak{L}_\hbar^2},
\end{aligned}$$
where $\CR^\hb[\ga;V]$ is defined in \eqref{eq: reformulated LS remainder}.
\end{lemma}

\begin{proof}
By H\"older's inequality, we have
$$\big|\langle \psi\big|(\nabla_q w \ast \zeta) f\big\rangle_{L_{q,p}^2}\big|\leq \|\psi\|_{L_{q,p}^2}\|\nabla_q w \ast \zeta\|_{L_q^\infty}\|f\|_{L_{q,p}^2}.$$
Similarly but with \eqref{eq: remainder estimate 1}, we obtain
$$\begin{aligned}
\big|\big\langle \psi \big|\CR^\hb[\ga;w \ast \zeta]\big\rangle_{L_{q,p}^2}\big|&\leq \|\langle\nabla_p\rangle \psi\|_{L_{q,p}^2}\big\|\langle \nabla_p\rangle^{-1}\mathcal{R}^\hbar[\gamma; w*\zeta]\big\|_{L_{q,p}^2(\mathbb{R}^6)}\\
&\lesssim \|\langle\nabla_p\rangle \psi\|_{L_{q,p}^2}\|\nabla_q w*\zeta\|_{L_q^\infty}\|\gamma\|_{\mathfrak{L}_\hbar^2}.
\end{aligned}$$
Hence, by the simple bound 
\begin{equation}\label{eq: basic derivative of potential bound}
\|\nabla_q w \ast \zeta\|_{L_q^\infty}\leq \|\nabla_q w\|_{L_q^\infty + L_q^{\frac{3}{1+a+\epsilon}}} \|\zeta\|_{L_q^1\cap L_q^{\frac{3}{2-a-\epsilon}}},
\end{equation}
the lemma follows.
\end{proof}



\begin{proof}[Proof of Proposition \ref{prop:Vlasov sol}]
\textbf{(Step 1. Weak solution on a compact interval)}
Fix any large finite $T>0$ and a sequence $\{\hbar_j\}_{j=1}^\infty$ such that $\hbar_j\to 0$.  We denote $f_j(t):=f^{\hbar_j}(t)$ and $g_j(t):=g^{\hbar_j}(t)=\mathcal{U}(-t)f^{\hbar_j}(t)$. By Proposition \ref{prop: compactness of fh}, passing to a subsequence, we may assume that $g_j(t)\to g(t)$ in $C_t([-T,T]; w-B_{\et})$, and let $f(t)=\mathcal{U}(t)g(t)$.

We will show that $g$ is a weak solution to \eqref{eq: Vlasov-moving frame'} in $C_t([-T,T]; H_{q,p}^{-1})$ in the sense that for any $\psi\in H_{q,p}^1$, 
\begin{equation}\label{eq: weak formulation of g equation}
\langle\psi|g(t)\rangle_{L_{q,p}^2}=\langle \psi|f(0)\rangle_{L_{q,p}^2}-\int_0^t \big\langle \nabla_p\CU(t_1)\psi\big| \nabla_q (w*\rho_f)f(t_1)\big\rangle_{L_{q,p}^2}dt_1
\end{equation}
holds in $C([-T,T])$. For the proof, we fix $\psi\in H_{q,p}^1$, and recalling the equation for $g_j(t)$, we write 
$$\begin{aligned}
\langle\psi|g_j(t)\rangle_{L_{q,p}^2}&=\langle\psi|f_j(0)\rangle_{L_{q,p}^2}-\int_0^t \big\langle \nabla_p\CU(t_1)\psi\big|\nabla_q (w*\rho_{f_j})f_j(t_1)\big\rangle_{L_{q,p}^2}dt_1\\
&\quad-i\int_0^t \big\langle\CU(t_1)\psi\big|\mathcal{R}^{\hbar_j}[\gamma^{\hbar _j}; w*\rho_{f_j}](t_1)\big\rangle_{L_{q,p}^2}dt_1
\end{aligned}$$
Indeed, we have
$$\langle\psi|g_j(t)\rangle_{L_{q,p}^2}\to\langle\psi|g(t)\rangle_{L_{q,p}^2}\textup{ in }C([-T,T])\quad\textup{and}\quad \langle\psi|f_j(0)\rangle_{L_{q,p}^2}\to \langle\psi|f(0)\rangle_{L_{q,p}^2},$$
because $g_j(t)\to g(t)$ in $C_t([-T,T]; w-B_{\eta_0})$ by Proposition \ref{prop: compactness of fh} while $f_j(0)\rightharpoonup f_0$ in $L_{q,p}^2$ by the construction of initial data $f_0$ (see Remark \ref{rmk:construction of f_0}). Therefore, it remains to show that for each $t_1\in[-T,T]$, 
$$\begin{aligned}
A_j(t_1)&:=\Big\langle \nabla_p\CU(t_1)\psi\Big|\Big(\nabla_q (w*\rho_{f_j})f_j-\nabla_q (w*\rho_f)f\Big)(t_1)\Big\rangle_{L_{q,p}^2}\to 0,\\
B_j(t_1)&:=\big\langle\CU(t_1)\psi\big|\mathcal{R}^{\hbar_j}[\gamma^{\hbar _j}; w*\rho_{f_j}](t_1)\big\rangle_{L_{q,p}^2}\to 0.
\end{aligned}$$
For these, we recall from Theorem \ref{thm: uniform bounds for nonlinear solutions} that 
$$\sup_{j\geq1}\sup_{t\in\mathbb{R}}\Big\{\|\rho_{f_j(t)}\|_{ L_{q}^1}+\langle t\rangle^{1+a+\epsilon}\|\rho_{f_j(t)}\|_{L_{q}^{\frac{3}{2-a-\epsilon}}}\Big\}\lesssim \eta_0,$$
because $\rho_{\gamma}^\hbar(x)=\rho_{\textbf{Wig}^\hbar[\gamma]}(x)$. On the other hand, by Proposition \ref{prop: compactness of fh}, we have $\|\rho_{f(t)}\|_{L_{q}^1} \lesssim \eta_0$ and 
$$\begin{aligned}
\|\rho_{f(t)}\|_{ L_{q}^{\frac{3}{2-a-\epsilon}}}&\leq \|\langle p\rangle^{-(1+a+2\epsilon)}\|_{L_{p}^{\frac{3}{1+a+\epsilon}}}\|\langle p\rangle^{1+a+2\epsilon} f(t)\|_{L_{q,p}^{\frac{3}{2-a-\epsilon}}}\\
&\lesssim_\epsilon\|\langle p\rangle^{1+a+2\epsilon} g(t)\|_{L_{q,p}^{\frac{3}{2-a-\epsilon}}}\lesssim \eta_0.
\end{aligned}$$
Hence, given any small $\delta>0$, employing Lemma \ref{lemma: integral estimates for derivation of VP}, we may replace $w$ and $\psi$ in $A_j(t_1)$ and $B_j(t_1)$ by $w_\delta\in C_c^\infty(\mathbb{R}^3)$ and  $\psi_\delta\in C_c^\infty(\mathbb{R}^6)$, respectively, up to $O(\delta)$-errors. Thus, it can be reduced to show that as $j\to\infty$,
$$\begin{aligned}
A_j^\delta(t_1)&:=\Big\langle \nabla_p\mathcal{U}(t_1)\psi_\delta\Big|\Big(\nabla_q w_\delta \ast \rh_{f_j}) f_j-(\nabla_q w_\delta \ast \rh_{f}) f\Big)(t_1)\Big\rangle_{L_{q,p}^2}\to0,\\
B_j^\delta(t_1)&:=\big\langle \mathcal{U}(t_1)\psi_\delta\big|\CR^{\hb_j}[\ga^{\hb_j};w_\delta \ast \rh_{f_j}](t_1)\big\rangle_{L_{q,p}^2}\to 0.
\end{aligned}$$
Indeed, by \eqref{eq: remainder estimate 2}, $B_j^\delta(t_1)\to 0$. It remains to consider $A_j^\delta(t_1)$.

Decomposing the difference of the nonlinear terms, we write 
$$\begin{aligned}
A_j^\delta(t_1)&=\Big\langle \nabla_p\mathcal{U}(t_1)\psi_\delta\Big|\Big(\nabla_q w_\delta \ast \rh_{(f_j-f)}\Big) f_j(t_1)\Big\rangle_{L_{q,p}^2}\\
&\quad+\Big\langle \mathcal{U}(-t_1)\Big(\nabla_q w_\delta \ast \rh_{f(t_1)}\cdot\nabla_p\mathcal{U}(t_1)\psi_\delta\Big)\Big|(g_j-g)(t_1)\Big\rangle_{L_{q,p}^2}\\
&=:A_{j;1}^\delta(t_1)+A_{j;2}^\delta(t_1).
\end{aligned}$$
Then, by the weak convergence of $g_j$ (Proposition \ref{prop: compactness of fh}), we see that $A_{j;2}^\delta(t_1)\to 0$. On the other hand, for $A_{j;1}^\delta(t_1)$, we note that $f_j$ is uniformly bounded in $L_{q,p}^2$ (see \eqref{fh L^2 norm}), $\nabla_p\mathcal{U}(t_1)\psi_\delta$ is compactly supported in $\mathbb{R}^6$ for each $t_1$, and $\|w_\delta \ast \rh_{f_j}\|_{C_q^k}, \|w_\delta \ast \rh_{f}\|_{C_q^k}\lesssim_{k,\delta}1$ for any $k\geq1$. Hence, the proof of  $A_{j,1}^\delta(t_1)\to 0$ can be reduced to that of the point-wise convergence 
\begin{equation}\label{eq: regularized pointwise force field convergence}
\lim_{j\to \infty}\nabla_q w_\delta \ast (\rh_{f_j}-\rho_f)(x)=0
\end{equation}
for every $x$. Indeed, by Lemma \ref{lemma: regularization of Weyl and Wigner transforms}, we observe that
$$\begin{aligned}
w_\delta*(\rho_{\mathbf{Hus}^{\hbar}[\gamma^{\hbar}]}-\rho_{f^{\hbar}})&=w_\delta*\Big(\frac{e^{-|\cdot|^2/\hbar}}{(\pi\hbar)^{3/2}}*\rho_{f^{\hbar}}-\rho_{f^{\hbar}}\Big)\\
&=\bigg(\frac{e^{-|\cdot|^2/\hbar}}{(\pi\hbar)^{3/2}}*w_\delta-w_\delta\bigg)*\rho_{f^{\hbar}}.
\end{aligned}$$
Thus, replacing the Wigner transform $f_j=\mathbf{Wig}^{\hbar_j}[\gamma^{\hbar_j}]$ in $\rho_{f_j}$ by the Husimi transform, \eqref{eq: regularized pointwise force field convergence} can be reduced further to show
\begin{equation}\label{eq: regularized pointwise force field convergence, Husimi version}
\lim_{j\to\infty}\nabla_q w_\delta*\Big(\rho_{\mathbf{Hus}^{\hbar_j}[\gamma^{\hbar_j}]}-\rho_f\Big)(x)=0
\end{equation}
for each $x\in\mathbb{R}^3$, since the difference can be estimated by \eqref{eq: basic derivative of potential bound}.

Now, we note that by H\"older's inequality with the momentum bounds (\eqref{eq: weighted bound for Husimi (1)} and \eqref{eq: weighted norm bound for the weak limit f}), one can show that 
\begin{equation}\label{eq: density estimate for C_1^h (2)}
\quad \|\rho_{(1-\chi_R(p))\mathbf{Hus}^{\hbar_j}[\gamma^{\hbar_j}](t)}\|_{L_q^{\frac{3}{2-a-\epsilon}}}, \|\rho_{(1-\chi_R(p))f(t)}\|_{L_q^{\frac{3}{2-a-\epsilon}}}\leq \delta,
\end{equation}
where $\chi_R(p)$ is a smooth cut-off  such that $\chi_R\equiv1$ if $|p|\leq R$ and $\chi_R\equiv0$ if $|p|\geq 2R$. Thus, by \eqref{eq: density estimate for C_1^h (2)} and the weak convergence $\mathbf{Hus}^{\hbar_j}[\gamma^{\hbar_j}]=G^\hbar*\mathcal{U}(t)g_j\rightharpoonup f=\mathcal{U}(t)g$ in $L_{q,p}^2$, we conclude that 
$$\begin{aligned}
&\lim_{j\to \infty}\nabla_q w_\delta*\Big(\rho_{\mathbf{Hus}^{\hbar_j}[\gamma^{\hbar_j}]}-\rho_f\Big)(x)\\
&=\lim_{j\to \infty}\nabla_q w_\delta \ast \rho_{\chi_R(p)(\mathbf{Hus}^{\hbar_j}[\gamma^{\hbar_j}]-f)}(x)+O(\delta)\\
&=\lim_{j\to \infty}\iint_{\mathbb{R}^6}\nabla_q w_\delta(x-q)\chi_R(p) \big(\mathbf{Hus}^{\hbar_j}[\gamma^{\hbar_j}]-f)(q,p)dqdp+O(\delta)=O(\delta).
\end{aligned}$$
Since $\delta>0$ is arbitrary, we obtain \eqref{eq: regularized pointwise force field convergence, Husimi version}, which in conclusion, completes the proof of \eqref{eq: weak formulation of g equation}.\\
\textbf{(Step 2. Upgrade to the global $H_{q,p}^{-1}$-solution)}
In \eqref{eq: Vlasov-moving frame'}, we note that $f(0)\in L_{q,p}^2$ and $g(t)\in C_t([-T,T]; w-B_{\eta_0})$. For the nonlinear integral term, we observe that by \eqref{eq: Vlasov-moving frame'}, \eqref{eq: weighted norm bound for the weak limit f} and \eqref{eq: decay of limit of densities}
\begin{equation}\label{eq: nonlinear term H^{-1}-bound}
\begin{aligned}
&\bigg\|\int_0^t(-t_1 \na_q + \na_p)\mathcal{U}(-t_1)\big(\nabla_q (w*\rho_f)\cdot f\big)(t_1) dt_1\bigg\|_{C([-T,T]; H_{q,p}^{-1})}\\
&\leq\int_{-T}^T \lg t_1\rg \| \nabla_q (w*\rho_f)(t_1) f(t_1)\big\|_{L^2_{q,p}} dt_1 \\
&\leq \int_{-T}^T \lg t_1 \rg\|\nabla_q (w*\rho_f)(t_1)\|_{L_q^\infty}\|f(t_1)\|_{L_{q,p}^2} dt_1
\lesssim \int_{-T}^T \frac{\eta_0^2}{\langle t_1\rangle^{a}} dt_1\lesssim \eta_0^2.
\end{aligned}
\end{equation}
Subsequently, we prove that the equation
$$g(t)=f_0+\int_0^t \CU(-t_1) \big(\nabla_q (w*\rho_f)\cdot\nabla_pf\big)(t_1)dt_1
$$
holds in $C([-T,T]; H_{q,p}^{-1})$. However, since \eqref{eq: nonlinear term H^{-1}-bound} confirms that all terms in the equation are bounded uniformly in $T$, we conclude that $g(t)$ solves the equation \eqref{eq: Vlasov-moving frame'} in $C(\mathbb{R}; H_{q,p}^{-1})$.\\
\textbf{(Step 3. Scattering)}
For convenience, we only consider the positive time direction. Let $g(t)$ be the solution in Proposition \ref{prop:Vlasov sol}. Then, by the equation \eqref{eq: Vlasov-moving frame'} and repeating the estimates in \eqref{eq: nonlinear term H^{-1}-bound}, we prove that as $t_2>t_1\to +\infty$,
$$\begin{aligned}
\|g(t_2)-g(t_1)\|_{H_{q,p}^{-1}}&=\bigg\|\int_{t_1}^{t_2} (-t'\na_q + \na_p)\CU(-t') \big(\nabla_q (w*\rho_f)\cdot f\big)(t')dt'\bigg\|_{H_{q,p}^{-1}} \\
&\lesssim\int_{t_1}^{t_2} \frac{\eta_0^2}{\langle t'\rangle^{a}} dt'\to0.
\end{aligned}$$
Therefore, we conclude that there exists $f_+\in H_{q,p}^{-1}$ such that $g(t) \to f_+$ in $H_{q,p}^{-1}$ as $t\to+\infty$.
\end{proof}

\subsection{Correspondence of the scattering states: Proof of Theorem \ref{thm: semiclassical limit} \texorpdfstring{$(iii)$}{(iii)}}\label{sec: correspondence of the scattering states}

Finally, we prove the correspondence of the scattering states between quantum and classical scattering states.

\begin{proposition}[Semi-classical limit of quantum scattering states]\label{prop:correspondence of the scattering state}
Using the notations in Proposition \ref{prop: compactness of fh}, let $\ga_\pm^\hb$ be the quantum scattering states and let $f_\pm$ be the classical scattering states given in Theorem \ref{thm: semiclassical limit}. Then, as $j\to\infty$, 
$$\textup{\textbf{Wig}}^{\hbar_j}[\ga^{\hb_j}_\pm]\rightharpoonup f_\pm\mbox{ in }L_{q,p}^2(\mathbb{R}^6).$$
\end{proposition}

For the proof, we need the following commutator estimate.
\begin{lemma}[Commutator estimate for the Weyl quantization]\label{lem:commutator estiamte}
If $f \in H^1_{q,p}(\R^6)$ and $V \in \dot{W}^{1,\infty}_x(\R^3)$, then
$$\big\|\big[V,\U \Weyl[f] \Uinv\big]\big\|_{\mathfrak{L}_\hbar^2} \le \hb \bt \|V\|_{\dot{W}^{1,\I}_x(\R^3)} \|f\|_{H^1_{q,p}(\R^6)}.$$
\end{lemma}

\begin{proof}
We denote $\ga(t):= \U \Weyl[f] \Uinv$. Observe that by the mean value theorem,
$$[V,\ga](x,x')=(V(x)-V(x'))\ga(x,x')=\nabla V(x_*)\cdot(x-x')\ga(x,x')=\nabla V(x_*)[x,\gamma](x,x')$$
for some $x_*$. Hence, it follows that 
$$\|[V,\gamma]\|_{\mathfrak{L}_\hbar^2} \le  \|\na V\|_{L^\I_x}\|[x,\gamma]\|_{\mathfrak{L}_\hbar^2}.$$
On the other hand, by Lemma \ref{lem:equation for f hbar} and the definition of the Weyl quantization, one can show that 
$$[x,\gamma]=\big[x,\Weyl[\E f]\big] = i\hbar \Weyl[\na_p \E f].$$
Hence, by the isometry of the Weyl quantization (Lemma \ref{lemma: basic properties of Wigner transform}), we prove that $\|[x,\ga]\|_{\FL^2_\hb} = \hb \|\na_p \CU(t)f\|_{L^{2}_{q,p}} \le \hb \bt \|f\|_{H^{1}_{q,p}}$.
\end{proof}

\begin{proof}[Proof of Proposition \ref{prop:correspondence of the scattering state}]
By symmetry, we only consider the positive time direction. Since by the isometry (Lemma \ref{lemma: basic properties of Wigner transform}), $\|\Wig[\ga_+^{\hb}]\|_{L_{q,p}^2}=\|\gamma_+^\hbar\|_{\mathfrak{L}_\hbar^2}=\|\gamma_0^\hbar\|_{\mathfrak{L}_\hbar^2}$ is bounded uniformly in $\hbar\in(0,1]$, it suffices to show 
\begin{equation}\label{eq: scattering state correspondence claim}
\lim_{j\to\infty}\big\langle  \ph \big| \textbf{Wig}^{\hbar_j}[\ga_+^{\hb_j}]-f_+\big\rangle_{L^2_{q,p}}\to0 
\end{equation}
for any $\varphi\in C_c^\infty(\mathbb{R}^6)$. For the proof of \eqref{eq: scattering state correspondence claim}, we fix any small $\delta>0$ and decompose 
\begin{equation*}
\begin{aligned}
\big\langle \ph \big| \textbf{Wig}^{\hbar_j}[\ga_+^{\hb_j}]-f_+\big\rangle_{L_{q,p}^2}
&=\big\langle\ph \big|\textbf{Wig}^{\hbar_j}\big[\ga_+^{\hb_j}-\mathcal{U}^{\hbar_j}(T)^*\ga^{\hb_j}(T)\mathcal{U}^{\hbar_j}(T)\big]
\big\rangle_{L_{q,p}^2}\\
&\quad + \big\langle\ph \big|
\textbf{Wig}^{\hbar_j}\big[\mathcal{U}^{\hbar_j}(T)^*\ga^{\hb_j}(T)\mathcal{U}^{\hbar_j}(T)\big] - \mathcal{U}(-T)f(T)
\big\rangle_{L^2_{q,p}}\\
&\quad + \big\langle\ph\big|\mathcal{U}(-T)f(T)-f_+\big\rangle_{L^2_{q,p}}\\
&=:\textup{A}_{T,j}+\textup{B}_{T,j}+\textup{C}_T,
\end{aligned}
\end{equation*}
where $T\gg1$ is a large number to be chosen later. 

For $\textup{A}_{T,j}$, by duality between the Wigner transform and the Weyl quantization (Remark \ref{remark: Wigner Weyl remark}), we write 
$$\textup{A}_{T,j}=(2\pi\hb_j)^3\Tr\Big(\textbf{Weyl}^{\hbar_j}[\ph]^*\Big(\ga_+^{\hb_j}-\mathcal{U}^{\hbar_j}(T)^*\ga^{\hb_j}(T)\mathcal{U}^{\hbar_j}(T)\Big)\Big).$$
Then, inserting the Duhamel formula for $\gamma^\hbar(t)$ and 
$$\gamma_+^{\hbar}=\gamma_0^{\hbar}-\frac{i}{\hbar}\int_0^\infty \mathcal{U}^{\hbar}(t_1)^*[\Phi^\hbar, \gamma^\hbar](t_1)\mathcal{U}^{\hbar}(t_1) dt_1,$$
we obtain 
$$\begin{aligned}
\textup{A}_{T,j}=-\frac{i}{\hbar_j}\int_T^\infty (2\pi\hbar_j)^3\textup{Tr}\Big(\mathcal{U}^{\hbar_j}(t_1)^*[\Phi^\hbar, \gamma^\hbar](t_1)\mathcal{U}^{\hbar_j}(t_1)\textbf{\textup{Weyl}}^{\hbar_j}[\varphi]\Big)dt_1.
\end{aligned}$$
In the above, by the cyclicity of the trace and Lemma \ref{lemma: Wigner transform of free evolution}, the trace term can be rearranged as 
$$\begin{aligned}
\textup{Tr}\Big(\mathcal{U}^\hbar(t_1)^*[\Phi^\hbar, \gamma^\hbar](t_1)\mathcal{U}^\hbar(t_1)\textbf{\textup{Weyl}}^{\hbar}[\varphi]\Big)&=\textup{Tr}\Big([\Phi^\hbar, \gamma^\hbar](t_1)\textbf{\textup{Weyl}}^{\hbar}[\mathcal{U}(t_1)\ph]\Big)\\ &=\textup{Tr}\Big(\gamma^\hbar(t_1)\Big[\textbf{\textup{Weyl}}^{\hbar}[\mathcal{U}(t_1)\ph], \Phi^\hbar(t_1)\Big]\Big).
\end{aligned}$$
Subsequently, we obtain
$$\begin{aligned}
|\textup{A}_{T,j}|\leq
\int_T^\infty\|\gamma^{\hbar_j}(t_1)\big\|_{\mathfrak{L}^2_\hb}
\frac{1}{\hbar_j}
\big\|\big[\textbf{Weyl}^{\hbar_j} [\mathcal{U}(t_1)\ph], \Phi^{\hbar_j}(t_1)\big]\big\|_{\mathfrak{L}^2_\hb}dt_1.
\end{aligned}$$
Note that by the conservation law, $\|\gamma^\hbar(t_1)\|_{\mathfrak{L}_\hbar^2}=\|\gamma_0^\hbar\|_{\mathfrak{L}_\hbar^2}\leq\eta_0$. On the other hand, Lemma \ref{lem:commutator estiamte} implies
$$\begin{aligned}
&\big\|\big[\textbf{\textup{Weyl}}^{\hbar}[\mathcal{U}(t_1)\ph], \Phi^\hbar(t_1)\big]\big\|_{\mathfrak{L}^2}
\le \hb \lg t_1 \rg \|\nabla\Phi^\hbar(t_1)\|_{L_x^\infty} \|\ph\|_{H^1_{q,p}}.
\end{aligned}$$
Therefore, collecting all and applying the uniform bound for $\Phi^\hbar$ (Theorem \ref{thm: uniform bounds for nonlinear solutions}), we prove that there exists a large $T\geq1$, independent of $j\geq1$,  such that  
$$|\textup{A}_{T,j}|\leq \eta_0\|\ph\|_{H^1_{q,p}}\int_T^\infty \frac{\eta_0}{\langle t_1\rangle^{a}}dt_1\leq\delta.$$
On the other hand, by the classical scattering (Proposition \ref{prop:Vlasov sol}), we have $|\textup{C}_T|\leq\delta$ for large $T\geq 1$, replacing $T$ by a larger number if necessary. 

It remains to consider $\textup{B}_{T,j}$ for fixed large $T\geq1$ as above. Indeed, by Lemma \ref{lemma: Wigner transform of free evolution}, we have
$$\textbf{Wig}^{\hbar}[\mathcal{U}^{\hbar}(T)^*\ga^{\hb}(T)\mathcal{U}^{\hbar}(T)]=\mathcal{U}(-T)\textbf{Wig}^{\hbar}[\ga^{\hb}(T)]=\mathcal{U}(-T)f^\hbar(T).$$
Then, it follows from the semi-classical limit result (Proposition \ref{prop: compactness of fh}) that for large $j$, 
$$|\textup{B}_{T,j}|=\big|\big\langle \ph \big| \mathcal{U}(-T)f^{\hbar_j}(T) - \mathcal{U}(-T)f(T)\big\rangle_{L^2_{q,p}}\big|=\big|\big\langle \ph \big| g^{\hbar_j}(T) -g(T)\big\rangle_{L^2_{q,p}}\big|\leq\delta.$$
Therefore, we conclude that $|\langle \ph|\textbf{Wig}^{\hbar_j}[\ga_+^{\hb_j}]-f_+\rangle_{L_{q,p}^2}|\leq 3\delta$. Since $\delta>0$ is arbitrary, we complete the proof.
\end{proof}

\subsection{A new scattering result of the Vlasov equation: Proof of Theorem \ref{thm:scattering Vlasov}} \label{sec:scattering Vlasov}
Define $\ga^\hb_0:=\Tp[f_0]$. Then, by Lemma \ref{lem:Schatten Toeplitz}, $\ga_0^\hb$ satisfies the assumptions in Theorem \ref{thm: semiclassical limit}.
Since $f_0^\hb:= \Wig[\ga_0^\hb] \rightharpoonup f_0$ in $L^2_{q,p}$ (see Remark \ref{remark: Wig Top} $(i)$), applying Theorem \ref{thm: semiclassical limit}, we obtain the result.

\section*{Acknowledgments}
We thank the anonymous referees for their careful review and constructive comments, which helped improve the quality and clarity of this manuscript. The first author was supported by JSPS KAKENHI Grant Number 24KJ1338. The second author was supported by National Research Foundation of Korea (NRF) grant funded by the Korean government (MSIT) (No. RS-2023-00219980, No. RS-2026-25479401).

\section*{Ethics declarations}
{\bf Conflict of interest} The authors declare that there is no conflict of interest.

\end{document}